\documentclass[12pt]{elsarticle}
\makeatletter
\def\ps@pprintTitle{%
 \let\@oddhead\@empty
 \let\@evenhead\@empty
 \def\@oddfoot{\centerline{\thepage}}%
 \let\@evenfoot\@oddfoot}
\makeatother

\usepackage{color} 
\usepackage{caption}
\usepackage{lscape}
\usepackage{afterpage}
\usepackage{pstricks}
\usepackage{pst-plot}
\usepackage{longtable}
\usepackage{dcolumn}
\usepackage{pst-node}
\usepackage{amsmath}
\usepackage{amssymb}
\usepackage{amsthm} 
\usepackage{multirow}
\usepackage{colortab}
\usepackage{color}
\usepackage{array}
\usepackage{colortbl}
\usepackage{textcomp}
\usepackage{pst-all}
\usepackage{subfigure}
\usepackage{tikz}
\usepackage{url}
\usepackage{hyperref}
\psset{arrows=->, labelsep=3pt, mnode=circle}

\usepackage{eurosym}

\usepackage[english]{babel}

\def\R{\mathbb{R}}

\newfont{\rams}{msbm10 scaled\magstep1}

\newcommand{\rio}{\mathbb{R}}

\begin{document}
\title{Robust sustainable development assessment with composite indices aggregating interacting dimensions: the hierarchical-SMAA-Choquet integral approach}

\author[Eco]{\rm Silvia Angilella}
\ead{angisil@unict.it}
\author[Eco]{\rm Pierluigi Catalfo}
\ead{pcatalfo.it}
\author[Eco]{\rm Salvatore Corrente}
\ead{salvatore.corrente@unict.it}
\author[Eco]{\rm Alfio Giarlotta}
\ead{giarlott@unict.it}
\author[Eco,por]{\rm Salvatore Greco}
\ead{salgreco@unict.it}
\author[Eco]{\rm Marcella Rizzo}
\ead{rizzom@unict.it}

\address[Eco]{Department of Economics and Business, University of Catania, Corso Italia, 55, 95129  Catania, Italy}
\address[por]{University of Portsmouth, Portsmouth Business School, Centre of Operations Research and Logistics (CORL), Richmond Building, Portland Street, Portsmouth PO1 3DE, United Kingdom}


\maketitle

\addcontentsline{toc}{section}{Abstract}

\begin{center}
{\large {\bf Abstract}} 
\end{center}

\noindent The evaluation of sustainable development -- and, in particular, rural development -- through composite indices requires taking into account a plurality of indicators, which are related to economic, social and environmental aspects. The points of view evaluated by these indices are naturally interacting: thus, a bonus has to be recognized to units performing well on synergic criteria, whereas a penalization has to be assigned on redundant criteria. An additional difficulty of the modelization is the elicitation of the parameters for the composite indices, since they are typically affected by some imprecision. In most approaches, all these critical points are usually neglected, which in turn yields an unpleasant degree of approximation in the computation of indices. In this paper we propose a methodology that allows one to simultaneously handle these delicate issues. Specifically, to take into account synergy and redundancy between criteria, we suitably aggregate indicators by means of the Choquet integral. Further, to obtain recommendations that take into account the space of fluctuation related to imprecision in nonadditive weights (capacity of the Choquet integral), we adopt the Robust Ordinal Regression (ROR) and the Stochastic Multicriteria Acceptability Analysis (SMAA). Finally, to study sustainability not only at a comprehensive level (taking into account all criteria) but also at a local level (separately taking into account economic, social and environmental aspects), we  apply the Multiple Criteria Hierarchy Process (MCHP). We illustrate the advantages of our approach in a concrete example, in which we measure the rural sustainability of 51 municipalities in the province of Catania, the largest city of the East Coast of Sicily (Italy).

\medskip 
 
\noindent{\bf Keywords}: {Sustainable development; composite indices; Choquet integral preference model;  Stochastic Multicriteria Acceptability Analysis; Robust Ordinal Regression; Necessary and Possible Preference.}

\section{Introduction}\label{introduction}
Sustainable development is a theme that is attracting more and more interest in experts, policymakers and laymen (see, e.g., \cite{atkinson2014handbook}). The starting point of the sustainable development discussion is usually fixed in 1987 with the Brundtland Report \cite{wced1987our}, according to which it concerns ``development that meets the needs of the present without compromising the ability of future generations to meet their own needs". However, the concept of sustainable development has a quite rich list of predecessors. In 1972, the  publication \textit{The Limits to Growth} \cite{meadows1972limits} used for the first time the term ``sustainable'' in a modern sense, pointing out the need of ``a world system that is: (i) sustainable without sudden and uncontrolled collapse; (ii) capable of satisfying the basic material requirements of all of its people''. Going further back in time at the beginning of the eighteenth century, in his \textit{Sylvicultura Oeconomica} \cite{Sylviculturaoeconomica}, von Carlowitz  invited to ``act with nature, and not against it'', stating the principle that the ``conservation and cultivation of timber should be conducted so as to provide a continuous, persistent and sustaining utilization''. Also many economists have contributed to the current debate on sustainable development. Just to give some examples, let us mention Fran\c{c}ois Quesnay and the physiocrats \cite{quesnay1958franccois} with their idea that any wealth comes from the Earth; Malthus \cite{malthus1809essay}, who studied population's pressure on limited natural resources; Ricardo \cite{ricardo1815essay}, who postulated diminishing returns from land use; Pigou \cite{pigou2013economics}, who investigated concepts of economic externalities and proposed a taxation to correct the related inefficiency; Lotka and Volterra \cite{lotka1926elements,Variazioniefluttazioni}, who studied the limits of ecological systems through the dynamics of prey-predator models; Georgescu-Roegen \cite{EntropyLaw}, who applied the concept of entropy to economics on the basis of the consideration that all natural resources are irreversibly degraded as soon as used in economic activities.

From a political and operational point of view, after being proposed in the Brundtland Report, the concept of sustainable development has been reconsidered in several Summits: the 1992 United Nations Conference on Environment and Development in Rio de Janeiro, the World Summit on Social Development in Copenhagen in 1995, the 2002 World Summit on Sustainable Development in Johannesburg,  the  2012 United Nations Conference on Sustainable Development, again in Rio de Janiero. The sustainable development has been also the subject of the work of several projects of  national and international organisations and governments, such as the European Commission's ``Beyond GDP", the OECD's ``Measuring the Progress of Societies", the  Commission on the Measurement of Economic Performance and Social Progress (CMEPSP), generally referred to as the Stiglitz-Sen-Fitoussi Commission. Recently, on September 25th, 2015, the UN General Assembly voted the 17 Sustainable Development Goals (SDGs)  of the 2030 Development Agenda ``Transforming our world: the 2030 Agenda for Sustainable Development" \cite{un2015transforming}. The SDGs have been developed on the basis of the previous Millennium Development Goals (MDGs)  stated in Millennium Declaration adopted by the General Assembly of United Nations on September 8th, 2000 \cite{assembly2000united}. As the MDGs were organized in  eight goals with 21 targets and 48 indicators, SDGs associates 169 targets and 232 indicators to the 17 goals. All in all, it is apparent that an increasing attention is being devoted towards data monitoring different aspects of sustainable development.   

In fact, beyond the space of the official documents of the United Nations, the idea that indicators are of fundamental importance  has been clearly stated in the debate on  sustainable development since some time ago. For example, \cite{meadows1998indicators} claims that we need many indicators organized hierarchically in an information system and properly integrated to translate ultimate means (solar energy, the biosphere, earth materials, biogeochemical cycles) into ultimate ends (happiness, harmony, identity, fulfillment, self-respect, community) \cite{daly1973steady}. In an analogous perspective, \cite{cash2003knowledge} concludes that the contribution of research, innovation, monitoring, and assessment to sustainability should be directed to manage boundaries between knowledge and action, and has to be coordinated in an integrated knowledge systems that supplies pieces of information perceived as salient, legitimate and credible. Similar concerns are taken into account in the postnormal science approach \cite{funtowicz1990uncertainty}, in which the socio-environmental questions are seen as typical post-normal problems characterized by high decision stakes and high system uncertainty. These issues have to be handled taking into consideration data whose imprecision can be described by five aspects: numeral, unit, spread, assessment (which expresses salient qualitative information about the data), and pedigree (which represents an evaluative description of the mode the information is produced). The importance of assessing sustainable development by analysing many indicators with most advanced machine learning and data mining methodologies is also advocated in \cite{perez2014classification,shaheen2011mining}.  

From a more applicative point of view, there have been several contributions providing discussions and surveys on the use of data, information and indicators for sustainable development (see, for example, \cite{dahl2012achievements,mayer2008strengths,ness2007categorising,Singh:2012}). In \cite{costanza2016modelling} some very interesting reflections on the construction and the use of sustainable development indicators have been proposed. The basic point is that a collection of many diversified indicators to take under control, the so called dashboard approach, is essential but not sufficient. It is instead necessary to combine the dashboard of indicators with a meaningful aggregation in composite indices permitting to measure progresses towards the desired goals. In this perspective, construction of composite indices becomes of fundamental importance for sustainable development. In fact, composite indices are more and more adopted in many research and applicative areas ranging from economic development \cite{bandura2008survey} to innovation \cite{grupp2004indicators} or tourism \cite{mendola2017building}, but it is in the domain of sustainable development that they have got the main interest \cite{bohringer2007measuring,Singh:2012}. Several critical issues have to be taken into account in using composite indices, in general \cite{Greco_et_al_composite,nardo2005handbook}, and for sustainable development, in particular:
\begin{description}
\item[(1)] \textit{weighting:} which weights should be assigned to the single elementary indicators? How to elicit all the parameters necessary to apply the preference model on the basis of the composite index?
\item[(2)] \textit{aggregation of the elementary indicators:} should we use the classical weighted sum or some more complex aggregation procedure? 
\item[(3)] \textit{robustness:} how much stable are the final results in changing the weights of single indicators and other parameters of the preference model?
\item[(4)] \textit{structure:} hierarchical organization of the elementary indicators in dimensions, goal, targets and the like;
\item[(5)] \textit{participation:} involvement of experts and stakeholders in the construction of the composite indices.
\end{description}

With respect to the above points, there is a widespread consensus that concepts and methodologies developed by MCDA (Multiple Criteria Decision Analysis; for a state-of-the-art see \cite{GreFigEhr}) can give a useful support for applications of composite indices in sustainable development (\cite{cinelli2014analysis,greco2017multi,Munda2016,rowley2012aggregating}). In this perspective, we aim to propose the application of an advanced MCDA methodology, the hierarchical-SMAA-Choquet integral approach \cite{angilella2015robust}, to construct composite indexes taking into account the previous points in the following way. 

First of all, regarding point \textbf{(2)}, we shall employ an aggregation procedure that is slightly more complex than the usual weighted sum, namely, the \textit{Choquet integral} \cite{Choquet} (see \cite{Grabisch1996} and \cite{pinar2014constructing} for the application of the Choquet integral to MCDA and composite indices, respectively). The reason for our choice is that the weighted sum is not able to represent synergy and redundancy between elementary indicators, which instead are quite relevant in sustainable development. Indeed, suppose to consider two indicators representing low emissions of two specific pollutants, and that one pollutant exacerbates the bad effects of the others. This is a case of synergy between the two indicators and it is natural to give a bonus if there are low emissions for both the pollutants. Analogously, one could imagine other two indicators, the first one related again to low emissions of a pollutant and the second one measuring the funds devoted to research on mitigation of contamination from the same pollutant. In this case, there is a redundancy between the two indicators because, of course, one can expect that the research results in a decrease of the emissions for the considered pollutant. Therefore, in case of low emission and consistent fund for the research, in order to avoid an over-evaluation, it is reasonable reducing the score that would be given by the usual weighted sum. The Choquet integral permits to take into account synergy and redundancy between criteria called, in general, interactions between criteria. This is possible assigning a weight to each subset of indicators rather than to each single indicator by means of what technically is called a capacity. After this, the indicators are therefore aggregated, so that the Choquet integral can be considered an extension of the weighted sum permitting to take into account the interaction between indicators. 

With respect  to points \textbf{(1)} and \textbf{(5)}, in general, the weights considered by the Choquet integral can be elicited with the participation of the stakeholders, the policymakers and the experts. Indeed, the non additive ordinal regression \cite{angilella2004assessing,marichal2000determination} applies the basic idea of the ordinal regression  \cite{jacquet1982assessing} to the Choquet integral permitting to elicit the weights from some preference information supplied by the stakeholders. This can be expressed in terms of judgments such as unit $a$ is better than (indifferent to) unit $b$, or indicator $i_1$ is more important than (as important as) indicator $i_2$ or, also, there is a synergy (a redundancy) between indicators $i_3$ and $i_4$. 
 
With respect to point \textbf{(3)}, \textit{Robust  Ordinal Regression} (ROR, \cite{greco2008ordinal}) and \textit{Stochastic  Multicriteria  Acceptability  Analysis} (SMAA, \cite{Lahdelma,Lahdelma_S2}) can be applied to the Choquet integral (see \cite{angilella2010non} and \cite{Angilella2012}, respectively). Indeed, point \textbf{(3)}, refers to the variability of the ranking supplied by the composite index due to the plurality of weights and, more in general, parameters that can be compatible with the information supplied by the stakeholders, the policy makers and the experts. To handle these aspects of the procedure, ROR defines the necessary preference relation, which holds for all the sets of compatible weights (and, more in general, parameters), and the possible preference relation, which holds for at least one set of compatible weights: the result is a \textit{necessary and possible preference} \cite{GiaGre2013}. Instead, SMAA proposes to consider a probabilistic ranking represented in terms of Rank Acceptability Indices (RAIs) and Pairwise Comparison Indices (PCIs). A RAI gives the probability that, picking randomly one set of compatible weights, a given unit attains a certain rank position, while PCI gives the probability that one unit gets a better value than another. The application of the SMAA methodology to composite indices has been proposed in \cite{greco2017stochastic}. 

Finally, with respect to  point \textbf{(4)}, the idea of a hierarchy of indicators is very rooted in the discussion about sustainable development since its origins. Indeed, following the  paradigm proposed by Passet -- sustainable development as intersection and relationship between the spheres of environment, economy and society \cite{passet1979economique} -- very usually, indicators are subdivided between these three domains. The above mentioned MDGs and SDGs are also organized in a hierarchy of goals, targets and indicators. \cite{costanza2016modelling}  proposes to harmonize the three spheres of economic, social and environmental aspects with the SDGs, by assigning the 17 goals as means of the higher level objectives of ``Efficient allocation'', ``Fair distribution'' and  ``Sustainable scale''. In MCDA the idea of a hierarchy of criteria has been discussed in \cite{CGShierarchy,CORRENTE20161}, where the authors propose an approach called the \textit{Multiple  Criteria  Hierarchy  Process} (MCHP). Applying MCHP to the Choquet integral coupled with ROR and SMAA, one can get therefore the hierarchical-SMAA-Choquet integral Approach \cite{angilella2015robust} that we are proposing for composite indices, in general, and for sustainable development, in particular. This will permit to simultaneously take into account all the above points \textbf{(1)}-\textbf{(5)}, which, as shown, have been already considered singularly in the literature and in the practice of composite indices, but have not been considered together yet.

The paper is organized as follows.  In Section 2 we deal with methodological aspects, introducing all basic notions relative to the Choquet integral preference model, MCHP, ROR, and SMAA. In Section 3, to show how the proposed methodology works in a real world problem, we apply it to a case study in which 51 municipalities within the province of Catania (Italy) are compared and ranked according to their level of rural sustainable development. Conclusions and future directions of research are summarized in Section 4. 

\section{Methodologies}
This preliminary section collects all the technical notions and methodologies that are employed in the paper. Specifically, after a brief introduction on terminology and notation (Section~\ref{preliminaries}), we recall the preference model based on the Choquet integral (Section~\ref{Choquet}), the robust ordinal regression applied to the Choquet integral (Section~\ref{NAROR_des}), the stochastic multicriteria acceptability analysis (Section~\ref{SMAA_des}),  and the hierarchical Choquet integral preference model (Section~\ref{MCHP_Ch}). These methodologies are then put together in Section~\ref{applied method 1}, which describes the approach used for the evaluation of the rural sustainability of municipalities. 

\subsection{Preliminaries} \label{preliminaries}
Let $A$ be a finite set of alternatives, and $G=\left\{g_1,\ldots,g_n\right\}$ a family of evaluation criteria, that is, maps $g_i \colon A \to \R$ for all $i \in I = \{1,\ldots,n\}$. Of course, in case of composite indicators, criteria have to be indentified with elementary indicators. We denote by $2^{G}$ the family of all subsets of $G$ (the \textit{powerset} of $G$). In this paper we assume, without loss of generality, that all criteria are of ``gain type": thus, the larger the evaluation, the better the alternative.

The only ``objective" information stemming for the evaluation of the alternatives is the \textit{dominance relation} $D$ on $A$, denoted by $a D b$ if $a$ is at least as good as $b$ for all criteria, and $a$ is strictly better than $b$ for at least one criterion. Regrettably, the (good) feature of objectiveness that dominance possesses is balanced out by the (bad) feature of its poorness, since many alternatives are deemed incomparable according to it. The latter fact in turn prevents one from obtaining an ``effective recommendation" for the problem at hand. Therefore, an aggregation of all evaluations on criteria is necessary to get a more refined -- hence more informative -- relation on the set of alternatives. To this aim, the literature considers three different families of aggregation methods: (i) value functions~\cite{Keeney76}, (ii) outranking relations~\cite{roy96}, and (iii) decisions rules~\cite{greco2001rough}.

A \textit{value function} is a map $U \colon A \to [0,1]$, which assigns a real value to each alternative $a \in A$: this number provides an exact estimation of how ``good" the alternative $a$ is. The output of this approach is a \textit{total ranking} on $A$.\footnote{Technically, a value function may generate a \textit{total preorder} on $A$, which is a reflexive, transitive, and complete binary relation $R$ on $A$. (Recall that \textit{reflexive} means $a R a$ for all $a \in A$, \textit{transitive} means that $a R b$ and $b R c$ implies $a R c$ for all $a,b,c \in A$, and \textit{complete} means that $a R b$ or $b R a$ holds for all distinct $a,b \in A$.) The difference between a total ranking and a total preorder is the possibility of ex-equo, which is allowed in the latter and not in the former. However, probabilist reasons typically yield a total ranking on the set of alternatives.} 

An \textit{outranking relation} is a binary relation $S$ on the set $A$ of alternatives, where $aSb$ means that alternative $a$ is ``at least as good as" alternative $b$: roughly speaking, there are arguments supporting the fact that $a$ outperforms $b$, as well as no strong arguments against this conclusion. It follows that the output of the aggregation process using an outranking relation is a \textit{partial ranking} on $A$, in which both \textit{ex-equo} and incomparability are possible.

Finally, \textit{decision rules} link the comprehensive evaluation of an alternative with its performances on criteria by means of statements of the type ``if ..., then ...".
The advantage of this methodology is that the output is easily understandable by the DM.

The approach employed in this paper to aggregate evaluations is of type (i), in fact it uses a value function in terms of the Choquet integral.

\subsection{The Choquet integral preference model} \label{Choquet}
A \textit{capacity} on the power set of the family $G$ of criteria is a function $\mu \colon 2^{G} \to [0,1]$ with the following properties: $\mu (\emptyset )=0$ \textit{(base constraint)}, $\mu (G) = 1$ \textit{(normalization constraint)}, and $\mu (R) \leq \mu (S)$ for all $R,S \in 2^G$ such that $R \subseteq S$ \textit{(monotonicity constraint)}. 

Given a set $A$ of alternatives and a capacity $\mu$ on $G$, the \textit{Choquet integral} \cite{Choquet} is a function $C_\mu \colon A \to \R_{0}^{+}$, which evaluates each alternative according to $\mu$. The formal definition of the Choquet integral of $a \in A$ is the following: 
$$
\displaystyle C_{\mu }(a)=\overset{n}{\underset{i=1}{\sum }}\left[g_{(i)}(a)-g_{(i-1)}\left( a\right)\right] \mu \left( N_{i}\right),
$$ 

\noindent where the subscripts ${(\cdot)}$ of the criteria stand for a permutation of $I = \{1,\ldots,n\}$  in a way that $0=g_{(0)}(a) \leq g_{(1)}\left( a\right)\leq\ldots\leq  g_{(n)}\left( a \right)$, and $N_{i}=\left\{ (i), \ldots ,(n)\right\} \subseteq I$ for each $i \in I$. 

It is well-known that a capacity $\mu$ on $G$ can be equivalently given using its \textit{M\"{o}bius representation}~\cite{Rota}, which is the (unique) function $m \colon 2^{G} \to \R$ such that the equality

\begin{equation}\label{MuMo}
\mu (S)\,= \,\displaystyle\sum_{R \subseteq S}\, m(R)
\end{equation}
\noindent holds for each $S \in 2^G$. Then, all constraints can be restated as follows: 
\begin{itemize}
    \item $m(\emptyset)=0\,$;
    \item $\displaystyle\sum_{T\subseteq G}m(T)=1\,$;
    \item $\displaystyle\sum_{T\subseteq S}m(T\cup\{g_i\})\geq 0$ for all $g_i\in G$ and $S\subseteq G\setminus\{g_i\}\,$.
\end{itemize}
A useful feature of the M\"{o}bius representation of $\mu$ is that the Choquet integral of $a \in A$ can be equivalently written as~\cite{Gilboa:1994}:
$$
\displaystyle C_{\mu }(a) \, =\, \sum_{T \subseteq G} \, m(T) \, \min_{g_i \in T} g_{i}(a).
$$

Criteria naturally interact with each other. Thus, the importance of a criterion $g_i \in G$ depends not only on its own relevance, but also on its contribution to all coalitions of criteria. The \textit{Shapley value} \cite{Shapley} of a criterion $g_i \in G$ takes into account these features in the following way: 
\begin{equation}\label{shapley_ind}
    \varphi\left(\{g_i\}\right)\,=\,\sum_{T\subseteq G\setminus\{g_i\}} \frac{\big(|G-T|-1 \big)! \: |T|!}{|G|!} \: \Big(\mu(T\cup\{g_i\})-\mu(T)\Big).
\end{equation}
\noindent Similarly to the Shapley value, we can also evaluate the relationship between two distinct criteria $g_i,g_j \in G$ computing their \textit{interaction index}~\cite{Murofushi1993}: 
\begin{equation}\label{interaction_ind}
    \varphi \left( \{ g_i,g_j\} \right) \, = \, \sum_{T\subseteq G\setminus\{g_i,g_j\}} \frac{\big(|G - T| - 2 \big)! \: |T|!}{\big(|G| - 1 \big)!}\: \Big(\mu(T\cup\{g_i,g_j\})-\mu(T\cup\{g_i\})-\mu(T\cup\{g_j\})+\mu(T) \Big).
\end{equation}
\noindent Again, M\"{o}bius representations yield a simplification of formulas (\ref{shapley_ind}) and (\ref{interaction_ind}) (see~\cite{grabisch2000equivalent}):
\begin{equation}\label{shapley_ind_mob}
    \varphi \left(\{g_i\}\right)=\sum_{g_i \in A\subseteq G}\frac{m(A)}{|A|}
\end{equation}
\noindent and 
\begin{equation}\label{interaction_ind_mo}
    \varphi \left(\{g_i,g_j\} \right) =\sum_{\{g_i,g_j\}\subseteq A \subseteq G} \frac{m(A)}{|A|-1} \,.
\end{equation}

A direct application of the preference model based on the Choquet integral appears hardly feasible, since it would require the elicitation of $2^{|G|}-2$ parameters, namely, all the values $\mu(T)$ for $\emptyset \subsetneq T \subsetneq G$ (since $\mu(\emptyset)=0$ and $\mu(G)=1$ are fixed, due to the definition of capacity). Computing such a huge number of parameters is almost impossible, even for a rather small set of criteria. 

The notion of \textit{$q$-additive} capacity, introduced in~\cite{Grabisch:1997}, makes the model better suited for applications, insofar as it requires the elicitation of definitively fewer parameters. Recall that a capacity is $q$\textit{-additive} if $m(T)=0$ for all $T \subseteq G$ such that $|T|>q$. Said differently, given a suitable integer $q$ such that $1 \leq q < n$, the $q$-capacity $\mu(T)$ of any family $T \subseteq G$ of criteria is obtained, according to formula (\ref{MuMo}), by summing up the M\"{o}bius value $m(R)$ of all its subfamilies $R \subseteq T$ having cardinality bounded by $q$. Luckily, in the majority of real world applications, a value of $q = 2$ appears to be unrestrictive. The obvious advantage of a $2$-additive capacity is that it involves the elicitation of a rather limited number of parameters, namely, $n+\binom{n}{2}$: a value $m(\{g_i\}) $ for each singleton $g_i \in G$, and a value $m(\left\{ g_i,g_j\right\} )$ for each unordered pair $\{g_i,g_j\} \subseteq G$. 

Using the M\"{o}bius representation (measure) $m$ of a $2$-additive capacity $\mu$, the three constraints (base, normalization, and monotonicity) assume the following form:

\begin{itemize}
    \item[(C.1)] \textit{(base)} $m(\emptyset) = 0\,$;
    \item[(C.2)] \textit{(normalization)} $\displaystyle\sum_{g_i \in G} m(\{g_i\}) + \sum_{\{ g_i,g_j\} \subseteq G} m(\{ g_i,g_j\}) =1\,$;
    \item[(C.3)] \textit{(monotonicity)} for all $g_i \in G$ and $\emptyset \neq T \subseteq G \setminus \{g_i\}$, 
		$
		\left\{
		\begin{array}{l}
		m(\{g_i\}) \geq 0\,\\[1mm]
		m(\{g_i\}) + \displaystyle\sum_{g_j \in T} m(\{g_i,g_j\}) \geq 0\,.
		\end{array}
		\right.         
		$
\end{itemize}

\noindent As a consequence, the Choquet integral of $ a \in A$ can be computed in a simpler way as follows:

\begin{eqnarray} \label{Choquet_Mobius}
    C_{\mu }(a) &=&\; \sum_{g_i \in G} m(\{g_i\}) \: g_{i}(a) \; + \sum_{\{ g_i, g_j\} \subseteq G} m(\{ g_i,g_j\}) \, \min  \{g_{i}(a), g_{j}( a) \}. 
\end{eqnarray}
In this context, the two equations (\ref{shapley_ind_mob}) and (\ref{interaction_ind_mo}) -- respectively expressing the Shapley value of a criterion and the interaction index of a pair of criteria -- can be further simplified as follows:
\begin{equation} \label{Shapley_Mobius}
    \varphi(\{ g_i\}) \: = \: m(\{ g_i\}) + \sum_{g_j \in G \setminus \{g_i\}}  \frac{m(\{ g_i,g_j\})}{2}\,,
\end{equation}
\noindent and
\begin{equation} \label{Murofushi_Mobius}
    \varphi(\{g_i,g_j\}) \: = \: m(\{ g_i,g_j\}) \, .  
\end{equation}

\subsection{Non Additive Robust Ordinal Regression (NAROR)}\label{NAROR_des}
NAROR (\textit{Non Additive Robust Ordinal Regression}, see \cite{angilella2010non}) belongs to the family of ROR methods (\textit{Robust Ordinal Regression}, see \cite{CGKSml,greco2008ordinal}). As it happens in all ROR methods, in NAROR the DM is asked to provide preference information related to a subset $A^*\subseteq A$ of \textit{reference alternatives}, which (s)he supposed to know quite well. This type of information is given at two (somehow complementary) levels of complexity by (i) comparing alternatives to each other, and (ii) comparing intensities of preferences to each other.

Specifically, for each $a,b,c,d \in A^* \subseteq A$, the DM may provide the following pieces of information:    
\begin{itemize}
    \item  $a$ is \textit{preferred}  to $b$, denoted by $a\succ b$ (which yields the constraint $C_{\mu}(a)\geq C_{\mu}(b)+\varepsilon$);
    \item  $a$ is  \textit{indifferent} to $b$, denoted by $a\sim b$ (which yields $C_{\mu}(a)=C_{\mu}(b)$);
    \item   $a$ is preferred to $b$ \textit{more than} $c$ is preferred to $d$, denoted by $(a,b)\succ^{*}(c,d)$ (which yields $C_{\mu}(a)-C_{\mu}(b)\geq C_{\mu}(c)-C_{\mu}(d)+\varepsilon$ and $C_{\mu}(c)\geq C_{\mu}(d)+\varepsilon$);
    \item  the intensity of preference of $a$ over $b$ is the \textit{same as} the intensity of preference of $c$ over $d$, denoted by $(a,b)\sim^*(c,d)$ (which yields $C_{\mu}(a)-C_{\mu}(b)= C_{\mu}(c)-C_{\mu}(d)$ and $C_{\mu}(c)\geq C_{\mu}(d)+\varepsilon$).
\end{itemize} 

Furthermore -- differently from other ROR methods -- in NAROR the DM may also provide some preference information on criteria $g_i,g_j,g_l,g_k \in G$, typically:
\begin{itemize}
    \item criterion $g_i$ is \textit{more important} than criterion $g_j$, denoted by $g_{i}\succ g_{j}$ (which yields the constraint $\varphi(\{g_i\})\geq\varphi(\{g_j\})+\varepsilon$);
    \item criteria $g_i$ and $g_j$ are \textit{indifferent}, denoted by $g_{i}\sim g_{j}$ (which yields $\varphi(\{g_i\})=\varphi(\{g_j\})$);
    \item criteria $g_i$ and $g_j$ are \textit{positively (negatively) interacting} (which yields $\varphi(\{g_i,g_j\}) \geq \varepsilon \;(\leq -\varepsilon)$);
    \item criterion $g_i$ is preferred to criterion $g_j$ \textit{more than} criterion $g_l$ is preferred to criterion $g_k$, denoted by $(g_{i},g_{j})\succ^* (g_{l},g_{k})$ (which yields $\varphi(\{g_i\})-\varphi(\{g_j\})\geq\varphi(\{g_l\})- \varphi(\{g_k\})+\varepsilon$ and $\varphi(\{g_l\})\geq \varphi(\{g_k\})+\varepsilon$); 
    \item the difference of importance between $g_i$ and $g_j$ is the \textit{same as} the difference of importance between $g_l$ and $g_k$, denoted by $(g_{i},g_{j})\sim^* (g_{l},g_{k})$  (which yields $\varphi(\{g_i\})-\varphi(\{g_j\})=\varphi(\{g_l\})- \varphi(\{g_k\})$).
\end{itemize}

Notice that, in the above constraints, $\varepsilon$ is an \textit{auxiliary variable} being a nonnegative small number used to convert strict inequalities into weak inequalities: for example, the weak inequality $C_{\mu}(a)\geq C_{\mu}(b)+\varepsilon$ codifies the strict inequality $C_{\mu}(a)>C_{\mu}(b)$. The output of NAROR consists of a pair $\left(\succsim^N,\succsim^P\right)$ of preference relations on the set $A$ of alternatives: this is called a \textit{necessary and possible preference} (\textit{NaP-preference}, see \cite{CGGMM2017,Giarlotta2014,Giarlotta2015,GiaGre2013,GiaWat2017b,GiaWat2017}).\footnote{Technically, a NaP-preference is a pair $\left(\succsim^N,\succsim^P\right)$ of binary relations on $A$ satisfying the following properties: \textit{(base condition)} $\succsim^N$ is a preorder; \textit{(extension)} $\succsim^P$ extends $\succsim^N$; \textit{(transitive coherence)} for all $a,b,c \in A$, $a \succsim^N b \succsim^P c$ implies $a \succsim^P c$, and $a \succsim^P b \succsim^N c$ implies $a \succsim^P c$; \textit{(mixed completeness)} for all $a,b \in A$, $a \succsim^N b$ or $b \succsim^P a$. 

Under the \textit{Axiom of Choice}, a NaP-preference can be characterized by the existence of a nonempty family $\mathcal{T}$ of total preorders on $A$ such that $\succsim^N = \bigcap \mathcal{T}$ and $\succsim^P = \bigcup \mathcal{T}$ (see \cite{GiaGre2013} for details).} 
The \textit{necessary preference} $\succsim ^N$ and the \textit{possible preference} $\succsim ^P$ are defined as follows for each $a,b \in A$:
$$
\begin{array}{llllllll}
        a \succsim ^N b &&\stackrel{\mathrm{def}}{\Longleftrightarrow} && C_{\mu}(a) \geq C_{\mu}(b) &\text{for \textit{all} compatible capacities}, \\
        a \succsim ^P b &&\stackrel{\mathrm{def}}{\Longleftrightarrow} && C_{\mu}(a) \geq C_{\mu}(b) &\text{for \textit{at least one} compatible capacity},    
\end{array}
$$ 
where a \textit{compatible capacity} is a set of  M\"{o}bius measures (satisfying base, normalization and monotonicity properties) for which the preference information provided by the DM is restored.

Denoted by ${\cal E}^{DM}$ the set of all the above constraints translating the DM's preference information (including constraints (C.1), (C.2), and (C.3), see Section~\ref{preliminaries}), the existence of a compatible capacity can be checked by solving the following linear programming problem:
$$
\varepsilon^*=\max\varepsilon \:, \qquad \mbox{subject to}\;\; {\cal E}^{DM}.
$$
If ${\cal E}^{DM}$ is feasible and $\varepsilon^{*}>0$, then there exists at least one compatible capacity. Otherwise, there are inconsistencies in the preferences provided by the DM: these inconsistencies can be identified (and removed) by means of one of the methods presented in \cite{mousseau2003resolving}.

In case that at least one compatible capacity exists, the necessary and possible preferences relative to each pair $a,b \in A$ of alternatives can be computed by using the following set of constraints: 

\begin{description}
    \item[${\cal E}^N(a,b):$] \quad $C_{\mu}(b)\geq C_{\mu}(a)+\varepsilon\;$ and $\;{\cal E}^{DM}$;
    \item[${\cal E}^P(a,b):$] \quad $C_{\mu}(a)\geq C_{\mu}(b)\;$ and $\;{\cal E}^{DM}$.
\end{description}

\noindent Specifically, $a$ is \textit{necessarily preferred} to $b$ ($a \succsim^N b$) if either ${\cal E}^{N}(a,b)$ is infeasible or $\varepsilon^{N}\leq 0$, where $\varepsilon^{N}=\max\varepsilon$, subject to ${\cal E}^{N}(a,b)$. Similarly, $a$ is \textit{possibly preferred} to $b$ ($a \succsim^P b$) if ${\cal E}^{P}(a,b)$ is feasible and $\varepsilon^{P}>0$, where $\varepsilon^{P}=\max\varepsilon$, subject to ${\cal E}^{P}(a,b)$.

\subsection{Stochastic Multiocriteria Acceptability Analysis (SMAA)}\label{SMAA_des}\label{SMAA}
As already summarized in Section~\ref{preliminaries}, the aggregation of the evaluations of alternatives on the considered criteria can be performed by either (i) a value function, or (ii) an outranking relation, or (iii) decision rules. In the first two approaches, several types of parameters need to be determined, namely: (1) a capacity on $G$ for the Choquet integral preference model; (2) indifference, preference and veto thresholds for the outranking approach~\cite{roy96}. 

Since in most cases the evaluations of the alternatives as well as the preference model parameters are not known with certainty, SMAA (\textit{Stochastic Multicriteria Acceptability Analysis}, see~\cite{Lahdelma,Lahdelma_S2}) methodology turn out to be extremely useful. In fact, SMAA is a family of MCDA methods, which take into account uncertainty and/or imprecision on both the preference parameters and the alternative evaluations: their goal is to provide a ``robust recommendation" for the solution of the problem at hand. 

Notice that SMAA methods have been already applied to choice, ranking and sorting problems by considering value functions and outranking relations as a preference model~\cite{tervonenfigueira}. In what follows, we briefly describe one of these methods, called SMAA-2~\cite{Lahdelma_S2}, which deals with ranking problems, and has a value function as its underlying preference model. We shall employ SMAA-2 in our approach. The rationale for this choice is related to the goal of \textit{ranking} the 51 municipalities in the province of Catania with respect to their environmental sustainability, since this method has features that fit quite well our setting. 

To make notation more compact, henceforth we shall denote alternatives by $a_{1},a_{2},\ldots$ (instead of $a,b, \ldots$ as before). The value function of each alternative is the weighted average of its numerical evaluations by all criteria. Thus, for each $a_k \in A$, we have 

\begin{equation}\label{SMAAfunction}
    U(a_k,w)\: = \: \sum_{i=1}^n w_{i} \, g_{i}(a_k)
\end{equation}
\noindent where $w\in W= \big\{(w_1,\ldots,w_n) \in \rio^n : \: w_i \geq 0 \text{ and } \sum_{i=1}^{n} w_{i}=1 \big\}$ is the vector of weights of criteria.

In SMAA methods, the preference information is represented by two probability distributions -- denoted by $f_{\chi}$ and $f_{W}$ -- defined, respectively, on the evaluation space $\chi=[g_{i}(a_k)]_{i,k}$ (which is composed of all the ``performance matrices" of alternatives) and on the weight space $W$. Then for each performance matrix $\xi\in\chi$ and each weight vector $w\in W$, equation (\ref{SMAAfunction}) yields a complete ranking of the alternatives.  As a consequence, the rank position of alternative $a_k \in A$ can be computed by the following \textit{rank function}: 
$$
rank(k,\xi,w)=1+\sum_{h\neq k}\rho\big(U(\xi_{h},w)>U(\xi_{k},w)\big),
$$
\noindent where $\rho(\mathrm{false})=0$ and $\rho(\mathrm{true})=1$. 

Given $\xi \in \chi$, SMAA-2 also computes the set of weights of criteria for which alternative $a_k$ assumes rank $s$. Formally, for each $s \in \{1,2,\ldots,|A|\}$, let
$$
W_{k}^{s}(\xi)=\left\{w\in W:rank(k,\xi,w)=s\right\}.
$$

In order to take into account the plurality of possible rankings and get robust recommendations, in SMAA-2 additional indices are computed, as described below.
\begin{itemize}
    \item \textit{Rank acceptability index $b^s(a_k)$}: this represents the probability that alternative $a_k$ has rank $s$. More precisely, $b^s(a_k)$ provides a measure of the set of parameters (compatible with the DM's preference information) that assigns rank $s$ to the alternative $a_k$. Formally, it is defined by
$$
b^{s}(a_k) \; = \; \int_{\xi\in \chi}f_{\chi}(\xi)\int_{w\in W_{k}^{s}(\xi)}f_{W}(w)\;dw\;d\xi \,.
$$
\item \textit{Pairwise winning index $p(a_h,a_k)$}: this represents the probability that alternative $a_h$ is preferred to alternative $a_k$. More precisely, $p(a_h,a_k)$ provides a measure of the set of parameters (compatible with the DM's preference information) for which $a_h$ is strictly better than $a_k$. Formally, it is defined by
$$
p(a_h,a_k) \; = \; \int_{w\in W} f_W(w) \: \int_{{\substack{\xi\in \chi:\;U(\xi_{h},w)> U(\xi_{k},w)}}} f_\chi(\xi) d\xi\; dw.
$$ 
\end{itemize}

From a computational point of view, the multidimensional integrals defining the above indices are estimated by using a Monte Carlo methodology.

\subsection{Multiple Criteria Hierarchy Process (MCHP) and the Choquet integral preference model}\label{MCHP_Ch}
\begin{center}
\begin{figure}
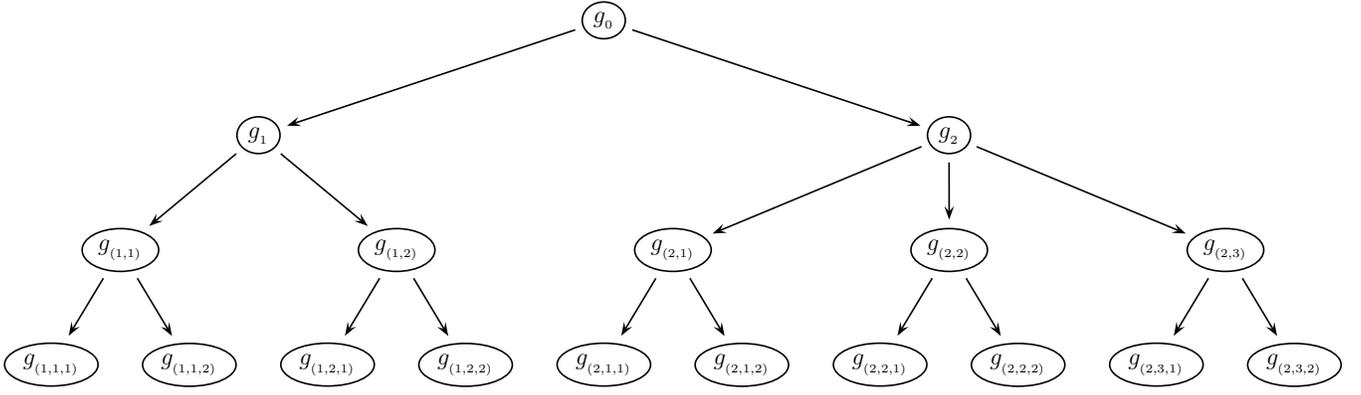

\resizebox{1\textwidth}{!}{
$
\pstree[nodesep=4pt,arrowsize=0.17]{\Toval{g_{_0}}}{%
\pstree{\Toval{g_{_{1}}}}{
\pstree{\Toval{g_{_{(1,1)}}}}{\Toval{g_{_{(1,1,1)}}} \Toval{g_{_{(1,1,2)}}}}
\pstree{\Toval{g_{_{(1,2)}}}}{\Toval{g_{_{(1,2,1)}}} \Toval{g_{_{(1,2,2)}}}}}
\pstree{\Toval{g_{_{2}}}}{%
\pstree{\Toval{g_{_{(2,1)}}}}{\Toval{g_{_{(2,1,1)}}} \Toval{g_{_{(2,1,2)}}}}
\pstree{\Toval{g_{_{(2,2)}}}}{\Toval{g_{_{(2,2,1)}}} \Toval{g_{_{(2,2,2)}}}}
\pstree{\Toval{g_{_{(2,3)}}}}{\Toval{g_{_{(2,3,1)}}} \Toval{g_{_{(2,3,2)}}}}
}
}
$
}
\caption{A hierarchy of criteria, displayed in three levels. {\small The root criterion is $g_{0}$ (at level zero), and has two subcriteria  $g_{1}$ and $g_{2}$ (at level one); there are ten elementary criteria at the last (third) level.\label{Hierarchy_ex}}} 
\end{figure}               
\end{center}

In MCHP (\textit{Multiple Criteria Hierarchy Process}, see~\cite{CGShierarchy}), the evaluation criteria are not all located at the same level, instead they are hierarchically structured. Thus, there is a \textit{root criterion} (the ``comprehensive objective") at level zero, a set of subcriteria of the root criterion at level one, etc.
The criteria at the lowest level of the hierarchy (which are the ``leaves" of the associated tree of criteria) are called \textit{elementary}. See Figure~\ref{Hierarchy_ex}.

To make terminology more suggestive, we shall be using the following notation: 
\begin{itemize} \small
    \item $G$ is the comprehensive set of criteria (at all levels of the hierarchy), and $g_0$ is the root criterion;
    \item $\textsf{I}_{G}$ is the set of indices of the criteria in $G$;
    \item $\textsf{E}_G\subseteq \textsf{I}_G$ is the set of indices of elementary criteria;
    \item $g_{\mathbf{r}}$ is a generic non-root criterion (where $\mathbf{r}$ is a vector with length equal to the level of the criterion);
    \item $g_{(\mathbf{r},1)},\ldots,g_{(\mathbf{r},n(\mathbf{r}))}$ are the immediate subcriteria of criterion $g_\mathbf{r}$ (located at the level below $g_{\mathbf{r}}$);
    \item $\textsf{E}(g_{\mathbf{r}})$ is the set of indices of all the elementary criteria descending from $g_{\mathbf{r}}$;
    \item $\textsf{E}(F)$ is the set of indices of the elementary criteria descending from at least one criterion in the subfamily $F \subseteq G$ (that is, $\textsf{E}(F) = \bigcup_{g_{\mathbf{r}}\in F} \textsf{E}(g_{\mathbf{r}})$);
    \item $G^{l}_{\mathbf{r}}$ is the set of subcriteria of $g_{\mathbf{r}}$ located at level $l$ in the hierarchy (below $g_{\mathbf{r}}$).
\end{itemize}
\noindent For instance, in Figure~\ref{Hierarchy_ex}, $\textsf{E}_G=\{(1,1,1),\ldots,(2,3,2)\}$, $g_{(\mathbf{1},1)}$ and $g_{(\mathbf{1},2)}$ are the immediate subcriteria of $g_{\mathbf{1}}$, $\textsf{E}(g_{\mathbf{2}})=\{(2,1,1),\ldots,(2,3,2)\}$, $G^{2}_{\mathbf{2}}=\{g_{(2,1)},g_{(2,2)},g_{(2,3)}\}$, $G^{3}_{\mathbf{2}}=\{g_{(2,1,1)},\ldots,g_{(2,3,2)}\}$, etc. 

Let $a \in A$ be an alternative, $g_{\mathbf{r}} \in G$ a non-elementary criterion (that is, $\mathbf{r}\in \textsf{I}_{G}\setminus \textsf{E}_G$), and $\mu$ a capacity defined on the powerset of of the set $\{g_{\mathbf{t}}: \mathbf{t}\in \textsf{E}_G\}$ of elementary criteria. The computation of the Choquet integral of $a$ on $g_{\mathbf{r}}$ is based on a capacity $\mu_{\mathbf{r}}$, which is derived from the original capacity $\mu$ and is defined on the powerset of the set of the immediate subcriteria of $g_{\mathbf{r}}$, that is, $\{g_{(\mathbf{r},1)},\ldots,g_{(\mathbf{r},n(\mathbf{r}))}\}$. Formally, the capacity $\mu_{\mathbf{r}}$ is defined as follows for each $F\subseteq\{g_{(\mathbf{r},1)},\ldots,g_{(\mathbf{r},n(\mathbf{r}))}\}$:
$$
\mu_{\mathbf{r}}(F)=\frac{\mu(\{g_{\mathbf{t}}: \mathbf{t}\in \textsf{E}(F)\})}{\mu(\{g_{\mathbf{t}}: \mathbf{t}\in \textsf{E}(g_{\mathbf{r}})\})}
$$
\noindent where $\;{\mu(\{g_\mathbf{t}: \mathbf{t}\in\textsf{E}(g_{\mathbf{r}})\},\emptyset)}\neq 0$ because, on the contrary, criteria from $\textsf{E}(g_{\mathbf{r}})$ would have null importance and would not be meaningful.

Then the Choquet integral of $a$ on $g_{\mathbf{r}}$ is computed by the formula
\begin{equation}
	C_{\mu_\textbf{r}}(a)=\frac{C_\mu(a_\textbf{r})}{\mu(\{g_{\mathbf{t}}: \mathbf{t}\in \textsf{E}(g_{\mathbf{r}}) \})}
	\label{ChoquetHierarchy}
\end{equation}

\noindent where $a_\textbf{r}$ is a ``fictitious" alternative having the same evaluations as $a$ on elementary criteria in $\textsf{E}(g_\textbf{r})$, and null evaluation on elementary criteria outside $\textsf{E}(g_\textbf{r})$ (that is, $g_\textbf{t}(a_\textbf{r})=g_\textbf{t}(a)$ if $\textbf{t}\in \textsf{E}(g_\textbf{r})$, and $g_\textbf{t}(a_\textbf{r})=0$ otherwise). 

As a consequence, we can naturally associate a preference relation $\succsim_{\mathbf{r}}$ to each non-elementary node $g_{\mathbf{r}}$ as follows:
$$
a\succsim_{\mathbf{r}}b \quad \stackrel{\mathrm{def}}{\Longleftrightarrow} \quad C_{\mu_\textbf{r}}(a)\geq C_{\mu_\textbf{r}}(b)
$$
\noindent where $a,b \in A$. 

\vspace{-1cm}
As previously observed, $2$-additive capacities are in general sufficient for the majority of practical cases. The Shapley value (\ref{Shapley_Mobius}) and the interaction index (\ref{Murofushi_Mobius}) -- which are classically defined for a ``flat" structure of criteria, that is, when no hierarchy is involved -- can be suitably computed in terms of the M\"{o}bius representation also for a hierarchical structure of criteria. For instance, given a non-elementary criterion $g_{\mathbf{r}}$ and three of its subcriteria placed at the level $l$ of the hierarchy -- say, $g_{(\mathbf{r},w)},g_{(\mathbf{r},w_1)},g_{(\mathbf{r},w_2)}\in G^{l}_{\mathbf{r}}$ -- we have: 

\begin{equation}\label{shapley_MCHP_2}
    \varphi_{\mathbf{r}}^{l}\left(\{g_{(\mathbf{r},w)}\}\right)
    = 
    \Bigg(\sum_{\mathbf{t}\in \textsf{E}(g_{(\mathbf{r},w)})}\!\!\!\!\!\!m\left(\{g_{\mathbf{t}}\}\right)
    +
    \!\!\!\!\!
    \sum_{\mathbf{t_1},\mathbf{t_2}\in \textsf{E}(g_{(\mathbf{r},w)})}\!\!\!\!\!\!\!\!\!m\left(\{g_{\mathbf{t_1}},g_{\mathbf{t_2}}\}\right)
    + 
    \!\!\!\!\!\!
    \sum_{\substack{\mathbf{t_{1}}\in \textsf{E}(g_{(\mathbf{r},w)}) \\ \mathbf{t_{2}}\in \textsf{E}({G}_{\mathbf{r}}^{l}\setminus\{g_{(\mathbf{r},w)}\})}}
    \!\!\!\!\!\!\!\!\!\!\!\!
    \frac{m(\{g_{\mathbf{t_{1}}},g_{\mathbf{t_{2}}}\})}{2}\Bigg)
    \frac{1}{\mu(\{g_{\mathbf{t}}: \mathbf{t}\in \textsf{E}(g_{\mathbf{r}}) \})},
\end{equation}

\noindent and
\begin{equation}\label{shapley_MCHP_3}
    \varphi_{\mathbf{r}}^{l}\left(\{g_{(\mathbf{r},w_1)},g_{(\mathbf{r},w_2)}\}\right) 
    \; =\; 
    \sum_{\substack{\mathbf{t_1}\in \textsf{E}\left(g_{(\mathbf{r},w_1)}\right)\\\mathbf{t_2}\in \textsf{E}\left(g_{(\mathbf{r},w_2)}\right)}} 
    \!\!\!\!
    m(\{g_{\mathbf{t_{1}}},g_{\mathbf{t_{2}}}\})\;
    \frac{1}{\mu(\{g_{\mathbf{t}}: \mathbf{t}\in \textsf{E}(g_{\mathbf{r}}) \})}\,.
\end{equation}

We refer the interested reader to ~\cite{angilella2015robust} for a detailed description of the extension of the Choquet integral preference model to the setting of hierarchically structured criteria.

\subsection{ROR and SMAA applied to the hierarchical Choquet integral preference model}
\label{applied method 1}
In order to apply the hierarchical Choquet integral preference model, we need to determine the M\"{o}bius representation of a capacity defined on the powerset of the set of elementary criteria. In particular, for the special case of a $2$-additive capacity, we have to determine $m(\left\{g_{\mathbf{t}}\right\})$ for each elementary criterion $g_{\mathbf{t}}$, and $m(\left\{g_{\mathbf{t_1}},g_{\mathbf{t_2}}\right\})$ for each unordered pair of elementary criteria $\left\{g_{\mathbf{t_1}},g_{\mathbf{t_2}}\right\}$. To elicit these values, we shall employ an ``ordinal regression" technique, which is based on the provision of suitable pieces of indirect preference information. Below we describe this technique.

Given a non-elementary criterion $g_{\mathbf{r}}$, the DM may provide the following type of preference information for some alternatives $a_h,a_k,a_v,a_z \in A$ (cf.\ Section~\ref{NAROR_des}):

\begin{itemize}
    \item $a_h$ is preferred to $a_k$ on $g_{\mathbf{r}}$, denoted by $a_h\succ_{\mathbf{r}}a_k$ (which is translated into the constraint $C_{\mu_{\mathbf{r}}}(a_h)\geq C_{\mu_{\mathbf{r}}}(a_k)+\varepsilon$);
    \item $a_h$ is indifferent to $a_k$ on $g_{\mathbf{r}}$, denoted by $a_h\sim_{\mathbf{r}}a_k$ (that is, $C_{\mu_{\mathbf{r}}}(a_h)=C_{\mu_{\mathbf{r}}}(a_k)$);
    \item on $g_{\mathbf{r}}$, $a_h$ is preferred to $a_k$ more than $a_v$ is preferred to $a_z$, denoted by $(a_h,a_k)\succ^{*}_{\mathbf{r}}(a_v,a_z)$ (that is, $C_{\mu_{\mathbf{r}}}(a_h)-C_{\mu_{\mathbf{r}}}(a_k)\geq C_{\mu_{\mathbf{r}}}(a_v)-C_{\mu_{\mathbf{r}}}(a_z)+\varepsilon\:$ and $\:C_{\mu_{\mathbf{r}}}(a_v)\geq C_{\mu_{\mathbf{r}}}(a_z)+\varepsilon$); 
    \item on $g_{\mathbf{r}}$, the intensity of preference of $a_h$ over $a_k$ is the same as the intensity of preference of $a_v$ over $a_z$, denoted by $(a_h,a_k)\sim_{\mathbf{r}}^{*}(a_v,a_z)$ (that is, $C_{\mu_{\mathbf{r}}}(a_h)-C_{\mu_{\mathbf{r}}}(a_k)=C_{\mu_{\mathbf{r}}}(a_v)-C_{\mu_{\mathbf{r}}}(a_z)$ and $C_{\mu_{\mathbf{r}}}(a_v)\geq C_{\mu_{\mathbf{r}}}(a_z)+\varepsilon$).
\end{itemize} 

Furthermore, considering criteria $g_{\mathbf{r_1}},g_{\mathbf{r_2}},g_{\mathbf{r_3}},g_{\mathbf{r_4}}\in G^{l}_{\mathbf{r}}$, the DM may provide the following type of preference information:

\begin{itemize}
    \item $g_{\mathbf{r_1}}$ is more important than $g_{\mathbf{r_2}}$, denoted by $g_{\mathbf{r_1}}\succ g_{\mathbf{r_2}}$ (which is translated into $\varphi_{\mathbf{r}}^{l}\left(\left\{g_{\mathbf{r_1}}\right\}\right)\geq\varphi_{\mathbf{r}}^{l}\left(\left\{g_{\mathbf{r_2}}\right\}\right)+\varepsilon$);
    \item $g_{\mathbf{r_1}}$ and $g_{\mathbf{r_2}}$ are equally important, denoted by $g_{\mathbf{r_1}}\sim g_{\mathbf{r_2}}$ (that is, $\varphi_{\mathbf{r}}^{l}\left(\left\{g_{\mathbf{r_1}}\right\}\right)=\varphi_{\mathbf{r}}^{l}\left(\left\{g_{\mathbf{r_2}}\right\}\right)$);
    \item $g_{\mathbf{r_1}}$ and $g_{\mathbf{r_2}}$ are positively interacting (that is, $\varphi_{\mathbf{r}}^{l}\left(\left\{g_{\mathbf{r_1}},g_{\mathbf{r_2}}\right\}\right)\geq\varepsilon$);
    \item $g_{\mathbf{r_1}}$ and $g_{\mathbf{r_2}}$ are negatively interacting (that is, $\varphi_{\mathbf{r}}^{l}\left(\left\{g_{\mathbf{r_1}},g_{\mathbf{r_2}}\right\}\right)\leq-\varepsilon$);
    \item the interaction between $g_{\mathbf{r_1}}$ and $g_{\mathbf{r_2}}$ is greater than the interaction between $g_{\mathbf{r_3}}$ and $g_{\mathbf{r_4}}$, which is codified differently according to the following subcases:
    \begin{itemize}
        \item if there is positive interaction between both pairs of criteria, then the constraints translating this preference are $\varphi_{\mathbf{r}}^{l}\left(\left\{g_{\mathbf{r_1}},g_{\mathbf{r_2}}\right\}\right)-\varphi_{\mathbf{r}}^{l}\left(\left\{g_{\mathbf{r_3}},g_{\mathbf{r_4}}\right\}\right)\geq\varepsilon$ and $\varphi_{\mathbf{r}}^{l}\left(\left\{g_{\mathbf{r_3}},g_{\mathbf{r_4}}\right\}\right)\geq\varepsilon$;
        \item  if there is negative interaction between both pairs of criteria, then the constraints translating this preference are $\varphi_{\mathbf{r}}^{l}\left(\left\{g_{\mathbf{r_1}},g_{\mathbf{r_2}}\right\}\right)-\varphi_{\mathbf{r}}^{l}\left(\left\{g_{\mathbf{r_3}},g_{\mathbf{r_4}}\right\}\right)\leq-\varepsilon$ and $\varphi_{\mathbf{r}}^{l}\left(\left\{g_{\mathbf{r_3}},g_{\mathbf{r_4}}\right\}\right)\leq-\varepsilon$;
    \end{itemize}
    \item $g_{\mathbf{r_1}}$ is preferred to $g_{\mathbf{r_2}}$ more than $g_{\mathbf{r_3}}$ is preferred to $g_{\mathbf{r_4}}$, denoted by $\left(g_{\mathbf{r_1}},g_{\mathbf{r_2}}\right)\succ^{*}\left(g_{\mathbf{r_3}},g_{\mathbf{r_4}}\right)$ (that is, $\varphi_{\mathbf{r}}^{l}(\{g_{\mathbf{r_1}}\}) -\varphi_{\mathbf{r}}^{l}(\{g_{\mathbf{r_2}}\}) \geq \varphi_{\mathbf{r}}^{l}(\{g_{\mathbf{r_3}}\}) - \varphi_{\mathbf{r}}^{l}(\{g_{\mathbf{r_4}}\}) +\varepsilon\:$ and $\:\varphi_{\mathbf{r}}^{l}(\{g_{\mathbf{r_3}}\}) - \varphi_{\mathbf{r}}^{l}(\{g_{\mathbf{r_4}}\}) \geq \varepsilon$);
    \item the difference of importance between $g_{\mathbf{r_1}}$ and $g_{\mathbf{r_2}}$ is the same of the difference of importance between $g_{\mathbf{r_3}}$ and $g_{\mathbf{r_4}}$, denoted by $\left(g_{\mathbf{r_1}},g_{\mathbf{r_2}}\right)\sim^{*}\left(g_{\mathbf{r_3}},g_{\mathbf{r_4}}\right)$ (that is, $\varphi_{\mathbf{r}}^{l}(\{g_{\mathbf{r_1}}\}) - \varphi_{\mathbf{r}}^{l}(\{g_{\mathbf{r_2}}\}) = \varphi_{\mathbf{r}}^{l}(\{g_{\mathbf{r_3}}\}) - \varphi_{\mathbf{r}}^{l}(\{g_{\mathbf{r_4}}\})$).
\end{itemize}

As explained in Section \ref{NAROR_des}, $\varepsilon$ is an auxiliary variable used to convert the strict inequalities into weak ones. Similarly, we denote by ${\cal E}^{DM}$ the set of constraints translating the DM's preference information along with the base, monotonicity and normalization constraints. All in all, to check if there exists at least one compatible capacity, one has to solve the following linear programming problem:
$$
\varepsilon^* \, = \, \max\varepsilon \,,\quad \mbox{subject to}\;\;{\cal E}^{DM}.
$$
\noindent If ${\cal E}^{DM}$ is feasible and $\varepsilon^{*}>0$, then there exists at least one compatible capacity; otherwise, there are inconsistencies, which need to be identified~\cite{mousseau2003resolving}. 

As in Section~\ref{NAROR_des}, the output of this analysis is a NaP-preference, whose computation is described below. Let $g_{\mathbf{r}}$ be a non-elementary criterion. Consider the following two sets of constraints: 
\begin{description}
    \item[${\cal E}^N_{\mathbf{r}} (a_h,a_k):$] \quad $C_{\mu_{\mathbf{r}}}(a_k)\geq C_{\mu_{\mathbf{r}}}(a_h)+\varepsilon\;$ and $\;{\cal E}^{DM}$;
    \item[${\cal E}^P_{\mathbf{r}} (a_h,a_k):$] \quad $C_{\mu_{\mathbf{r}}}(a_h)\geq C_{\mu_{\mathbf{r}}}(a_k) \;$ and $\;{\cal E}^{DM}$.
\end{description}

Then, $a_h$ is \textit{necessarily preferred} to $a_k$ with respect to criterion $g_{\mathbf{r}}$ (denoted by $a_h \succsim^N_{\mathbf{r}} a_k$) if either ${\cal E}_{\mathbf{r}}^{N}(a_h,a_k)$ is infeasible, or $\varepsilon_{\mathbf{r}}^{N}\leq 0$, where $\varepsilon_{\mathbf{r}}^{N}=\max\varepsilon$, subject to ${\cal E}_{\mathbf{r}}^{N}(a_h,a_k)$. Analogously, $a_h$ is \textit{possibly preferred} to $a_k$ with respect to criterion $g_{\mathbf{r}}$ (denoted by $a_h \succsim^P_{\mathbf{r}} a_k$) if both ${\cal E}_{\mathbf{r}}^{P}(a_h,a_k)$ is feasible and $\varepsilon_{\mathbf{r}}^{P}>0$, where $\varepsilon_{\mathbf{r}}^{P}=\max\varepsilon$, subject to ${\cal E}_{\mathbf{r}}^{P}(a_h,a_k)$.

In concrete scenarios, it is likely that simultaneously $a_h$ is possibly preferred to $a_k$, and $a_k$ is possibly preferred to $a_h$. However, in these cases it may also happen that the number of compatible capacities for which $a_h$ is preferred to $a_k$ is very different from the number of compatible capacities for which $a_k$ is preferred to $a_h$. In other words, the necessary and possible approach is invariant of the cardinality of compatible capacities. This is the point where the SMAA methodology turns to be useful, insofar as it allows one to estimate how good an alternative is in comparison to others, and how often it is preferred over the others. Below we explain how to apply SMAA to the hierarchical Choquet integral preference model. 

The set of linear constraints in ${\cal E}^{DM}$ defines a convex set of M\"{o}bius parameters. Then we may employ the \textit{Hit-And-Run} (HAR) method~\cite{smith1984,Tervonen2012,van_valkenhoef} to explore this set of parameters. Specifically, HAR iteratively samples a set of capacities expressed in terms of M\"{o}bius parameters satisfying ${\cal E}^{DM}$, until a stopping condition is met. For each sampled set of M\"{o}bius parameters and each given criterion $g_{\mathbf{r}}$, one can compute values of the Choquet integral for all alternatives. The output is a ranking of the alternatives with respect to criterion $g_{\mathbf{r}}$. Having many rankings as samples, one can then compute all indices that are typical of the SMAA methodology (see Section~\ref{NAROR_des}), namely: 

\begin{itemize}
    \item the \textit{rank acceptability index} $b^{s}_{\mathbf{r}}(a_k)$, which is the frequency of having alternative $a_k$ at position $s$ in the ranking obtained w.r.t.\ criterion $g_{\mathbf{r}}$;
    \item the\textit{ pairwise winning index} $p_{\mathbf{r}}(a_h,a_k)$, which is the frequency of having a preference of $a_h$ over $a_k$ w.r.t.\ criterion $g_{\mathbf{r}}$.
\end{itemize}

\noindent Furthermore, using the rank acceptability indices, two additional indices can be computed as follows~\cite{naror_TO,milos_tommi_1}:

\begin{itemize}
    \item the \textit{downward cumulative rank acceptability index} $b^{\leq s}_{\mathbf{r}}(a_k)$, which is the frequency of having alternative $a_k$ at a position not greater than $s$ w.r.t.\ criterion $g_{\mathbf{r}}$, that is,
    $$
    b^{\leq s}_{\mathbf{r}}(a_k)=\sum_{q=1}^{s}b^{q}_{\mathbf{r}}(a_k)\,;
    $$
    \item \textit{the upward cumulative rank acceptability index} $b^{\geq s}_{\mathbf{r}}(a_k)$, which is the frequency of having alternative $a_k$ at a position not lower than $s$ w.r.t.\ criterion $g_{\mathbf{r}}$, that is,
    $$
    b^{\geq s}_{\mathbf{r}}(a_k)=\sum_{q=s}^{n}b^{q}_{\mathbf{r}}(a_k)\,.
    $$
\end{itemize}
    
Finally, observe that the employed methodology allows us to get both NaP-preferences and SMAA indices at a comprehensive level as well (that is, with respect to the root criterion $g_{\mathbf{0}}$). 

\section{The case study}
To show how the hierarchical-SMAA-Choquet integral works, we apply the methodology described in the previous section to compare and rank with respect to their level of rural sustainable development 51 municipalities within the province of Catania, a city on the East Coast of the Italian region of Sicily. Before presenting, step by step, the application of the procedure in the case study, we give a brief sketch of the overall context in which the results of our analysis need to be interpreted, that is, the European Union (EU) policy with respect to the  rural sustainability development. After, we present the application of the proposed approach showing the support that it supplies to stakeholders, experts and policy makers. 

\subsection{The rural sustainable development in the European Union policy}
The structure of the Common Agricultural Policy (CAP), that is, the agricultural policy of the European Union (EU), is based on two different pillars. The first one is characterized by a specific \textit{sectoral policy}, in fact it is the pillar of the globalized agricultural market for commodities. The second one -- probably the most important and current -- is characterized by a \textit{territorial approach}, and it is useful to support the social cohesion and the integrated sustainable development of all rural areas~\cite{DGAgri}. With respect to this latter point, nowadays, it is apparent that the rural development policy is inspired by the strategic approach of EUROPA 2020 \cite{european2010europe}, which is structured to foster a smart, sustainable and inclusive development. As a matter of fact, sustainability has become the key factor in motivating all European policy actions aimed at sustaining any kind of human activity~\cite{buckwell2015should}. These actions are financially supported by both the EU and the national governments of its members, according to a co-financing principle. As a consequence, resources are distributed at a national and a regional level, according to specific criteria and priorities of intervention, and based on an estimation of the effectiveness of the employed policies~\cite{Nazzaro}. Consequently, rural development depends on all those policies that have as a goal the creation of a new balance between urban and territorial space, as well as the defence of the quality of rural life, including young farm generations and a long term sustainable environmental quality~\cite{rapisarda2014applicative}. Indeed, the conditions under which rural sustainability is defined, looks at a long time view through the logic of safeguarding and enhancing the resources of the territory, the environment and human capacities, including culture, the tradition of the territory and the complex of the intangible and relational resources \cite{Catalfo}. Therefore, the actual understanding of the concept of agriculture is determinate by the abandon of the vision of agriculture as a sector and has moved towards a territorial dimension characterized by a form of cross-sectoral rurality with reciprocal implications to other economic, social and environmental dimensions. In this context, the rural dimension looks at the role and the functions of agriculture. They emerge from the necessity to protecting the environment and the landscape, capitalizing on natural and environmental resources through tourism and culture, technological innovation, social cohesion and the development of quality of life in rural centers outside strictly urban areas \cite{PlataniaQuality}. 
 
This strategic policy framework also motivates the interest of the EU in controlling climate changes, mainly induced by a wrong productive system that badly affects the environment. In view of an effective financial strategy of support, it appears then inevitable to undertake a careful comparison of different rural areas with respect to several indicators that are naturally linked to sustainability \cite{teeter2016constraining}. Measuring sustainability has been already the goal of \textit{green accounting}, outperformed by evaluating several and different sustainability indexes~\cite{Munda2005}. However, there is not a unique system to combine such different types of indexes in a way  immediately  useful for policy makers. The sustainability indexes are typically of two categories: monetary and physical. For example, the  Index of Sustainable Economic Welfare (ISEW)\cite{Daly}, Weak Sustainability Index \cite{pearce_hamilton_atkinson_1996}, and the so-called El Serafy approach \cite{Yusuf}  belong to the first category of indexes, whereas  instances  of the second category are Human Appropriation of Net Primary Production \cite{Vitousek}, and the Ecological Footprint \cite{Rees}. All these indicators belong to the family of composite indicators, so that the approach we are proposing of a new generation methodology to construct composite indices seems well appropriate. 

\subsection{The application and the results of the hierarchical-SMAA-Choquet integral approach}
The specific example that we are going to analyze regards the evaluation of the rural sustainability of 51 municipalities within the province of Catania, which is the second largest city in the region of Sicily (Italy). To overcome the difficulties connected to the very definition of ``rural area", we have identified these municipalities using the classification in the 2007-2013 program of rural development of Sicily. The evaluation is based on a set of indicators that has been already used in~\cite{Boggia:2014} to assess the sustainable rural development of the municipalities in Umbria, another Italian region. In the latter case study, the evaluation has been performed using another methodology, namely, the Dominance-based Rough Set Approach (DRSA)~\cite{greco2001rough}. 

The reasons for which we choose a \textit{provincial} scale rather than a \textit{regional} scale -- which is the one typically employed in this type of analysis -- are connected to the possibility to directly interpret the results of the decision model in view of a deep knowledge and understanding of the territory. Technically, this change of scale is legitimated by the fact that the province of Catania has a level of complexity which is definitively comparable to that of a region, in terms of quality and variety of resources as well as on the level of social relationship. 

As confirmed by some previous research concerning the Eastern Coast of Sicily, the province of Catania offers a privileged position of observation for an assessment of rural sustainability, due to the contiguity of remarkably different socio-economic situations: residual but well-rooted industrial districts, tiny but highly developed technological areas, locations exclusively dedicated to touristic activity, agro-naturalistic settlements on the slopes of the Etna volcano, wineries that make the pride of the Sicilian island and so on. On the other hand, several areas of the province of Catania are characterized by a clear and marked backwardness, in particular rural ares with serious problems of underdevelopment, a fragile eco-environmental condition and a low level of anthropic dimension. \\
In view of all these highly contrasting and variegated features of the surroundings of the province of Catania, our analysis -- even if implemented at a provincial level -- has no methodological shortcomings.

Our original task was to apply the described methodology to evaluate the rural sustainable development of \textit{all} the 58 municipalities within the province of Catania (see OECD~\cite{OECD} and~EUROSTAT \cite{EU2001}). However, 7 out of these 58 municipalities were to be excluded from our analysis, since they are not classified as ``rural" according to the 2007-2013 program of the rural development of Sicily. The remaining 51 municipalities, considered in the presented study, are listed in Table~\ref{alternatives}.

\begin{table}[htbp]
\caption{List of the municipalities of Catania}\label{alternatives}\smallskip
\centering
\resizebox{0.8\textwidth}{!}{
\begin{tabular}{clclcl} 
\hline
\hline
\textbf{Alternative} & \textbf{Municipality} & \textbf{Alternative} & \textbf{Municipality} & \textbf{Alternative} & \textbf{Municipality} \\
\hline
\hline
$a_{1}$	& Aci Bonaccorsi         & $a_{18}$	& Linguaglossa & $a_{35}$	& Ramacca\\ 
$a_{2}$	& Aci Catena             & $a_{19}$	& Maletto      & $a_{36}$	& Randazzo\\ 
$a_{3}$ & Aci Sant'Antonio       & $a_{20}$	& Maniace      & $a_{37}$	& Riposto\\ 
$a_{4}$	& Acireale               & $a_{21}$	& Mascali      & $a_{38}$	& San Cono\\ 
$a_{5}$ &	Adrano                 & $a_{22}$	& Mascalucia   & $a_{39}$	& San Giovanni la Punta\\  
$a_{6}$	& Belpasso               & $a_{23}$	& Mazzarrone   & $a_{40}$	& San Michele di Ganzaria\\ 
$a_{7}$	& Biancavilla            & $a_{24}$	& Militello in Val di Catania & $a_{41}$	& San Pietro Clarenza\\ 
$a_{8}$	& Bronte                 & $a_{25}$	& Milo                        & $a_{42}$	& Santa Maria di Licodia\\ 
$a_{9}$	& Calatabiano            & $a_{26}$  & Mineo                      & $a_{43}$	& Santa Venerina\\
$a_{10}$	& Caltagirone          & $a_{27}$  & Mirabella Imbaccari        & $a_{44}$	& Sant'Alfio\\  
$a_{11}$	& Camporotondo Etneo   & $a_{28}$	 & Nicolosi                   & $a_{45}$	& Scordia\\ 
$a_{12}$	& Castel di Iudica     & $a_{29}$	 & Palagonia                  & $a_{46}$	& Trecastagni\\
$a_{13}$	& Castiglione di Sicilia & $a_{30}$	& Patern\'o                 & $a_{47}$	& Tremestieri Etneo\\ 
$a_{14}$	& Fiumefreddo di Sicilia & $a_{31}$	& Pedara                    & $a_{48}$	& Valverde\\ 
$a_{15}$	& Giarre                 & $a_{32}$	& Piedimonte Etneo          & $a_{49}$	& Viagrande\\ 
$a_{16}$	& Grammichele            & $a_{33}$	& Raddusa                   & $a_{50}$	& Vizzini\\
$a_{17}$	& Licodia Eubea          & $a_{34}$	& Ragalna                   & $a_{51}$	& Zafferana Etnea\\
\hline
\hline
\end{tabular}
}
\end{table}

The data set of the municipalities is collected from ISTAT (the Italian Institute of Statistics, \textit{www.istat.it}), considering the census data relative to 2010 (on agriculture) and 2011 (on population). The municipalities are evaluated on the basis of the criteria described in~\cite{Boggia:2014}. Specifically, three macro-criteria are considered: Social (So), Economic (Ec), and Environmental (En). These criteria have been further decomposed into more detailed sub-criteria as follows:

\begin{description}
\item[\textbf{Social Sustainability.}] In order to measure social sustainability in rural areas, three simple indicators were taken into consideration, since data used to compute them are very easy to be collected. 
    These indicators are useful to represent the presence of human settlements capable of maintaining an equilibrium between the original identity of the territory and its anthropic size. Notice that the large size of settlements has typically a negative impact on the rural dimension, due to an increased usage of resources and a consequent impoverishment of the territory~\cite{Gaviglio}.    
    \begin{enumerate}
    \item \textit{Population Scattering Index (PSI)} : it is the ratio between the population living in individual homes, villages, small towns, and the total population. An increase of this indicator signals an improvement of the social dimension for the sustainability of rural development. The relevance of PSI in the assessment of rural development lies in the fact that it is based on human settlements in terms of both territorial distribution and absorption of resources.  The presence of habitual settlements on a rural territory --which is by definition scarcely inhabited -- reveals the vitality of the rural area under consideration.
		\item \textit{Concentrated Population/km$^2$ (CP)}: This indicator measures the degree of housing density in the considered rural area. However, contrary to PSI, its increase highlights a worsening of the rural development, due to the motivations described for PSI. 
    \item \textit{Total Residents (R)}: Also for this indicator an inverse direction is observed, that is, the higher its level, and the lower the performance. Indeed, with respect to social sustainability, if the resident population increases, the level of rural development decreases, due to a depletion of the resources of the territory.
    \end{enumerate}
		
\item[\textbf{Economic Sustainability.}] Among the sustainability indicators of rural development related to the economic dimension, the following ones have been considered: 
    \begin{enumerate}
    \item \textit{Owned Homes/Total Homes (OH)}: According to the current literature \cite{Kaiser1990}, an increase of this index is correlated to a worsening of the rural development. Indeed, an enlargement of property ownership by residents is a condition of fragility for rural sustainability, because the growing number of people living in the area naturally induces the necessity of new houses. A value of the index close to one means that the trend to build new houses is booming. 
    \item \textit{Agricultural Age Structure (AA)}: This indicator is directly linked to the European strategy ``Europe 2020". An increase of this index represents a tendency for a greater entrepreneurial engagement capacity in a specific territory.
    \item \textit{Number of Bed Spaces in Rural Tourism Accommodation/km$^2$ (BS)}: With respect to this indicator, the areas with a higher number of farmhouses are more oriented to an  economic development.
    \end{enumerate}

\item[\textbf{Environmental Sustainability.}] Among several possible indicators, we have selected those that somehow take into account -- even in a hybrid way -- the fact that economy and environment  have to be considered together. Specifically, we consider:  
    \begin{enumerate}
    \item \textit{Typical Products Rate (TPR)}: It is ratio of  the number of typical products of the considered municipality over the total number of typical products in Sicily. A typical product  characterizes a municipality since its production is related to its specific environmental conditions ~\cite{PilatoQuality,PlataniaQuality2}.
    \item \textit{Irrigation Water Use/UAA (IWU)}.
    \item \textit{Livestock Standard units/UAA (LS)}.
    \item \textit{Organic UAA/UAA (OU)}. 
    \end{enumerate}

The classical indicators 2, 3 and 4 are related to the resource absorption of the surface devoted to agriculture~\cite{PetinoQuality}. As a measurement of this territorial dimension, we employed the Utilized Agricultural Area  (UAA), according to what is typical for this purpose.
\end{description}

\begin{figure}[!htb]
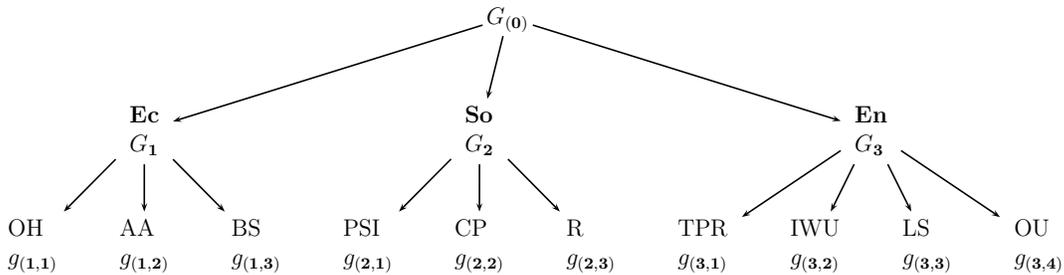

\centering
\caption{Hierarchical structure of criteria considered in the case study\label{hierarchy_sicily}}
\resizebox{0.8\textwidth}{!}{
\pstree[nodesep=2pt]{\TR{$G_{\mathbf{(0)}}$}}{       
				\pstree{\TR{$\begin{array}{l}\textbf{Ec}\\ G_{\mathbf{1}}\end{array}$}}{
        \TR{$\begin{array}{l}\text{OH}\\ g_{\mathbf{(1,1)}}\end{array}$}
				\TR{$\begin{array}{l}\text{AA}\\ g_{\mathbf{(1,2)}}\end{array}$}
				\TR{$\begin{array}{l}\text{BS}\\ g_{\mathbf{(1,3)}}\end{array}$}
				}
        \pstree{ \TR{$\begin{array}{l}\textbf{So}\\ G_{\mathbf{2}}\end{array}$}}{
				\TR{$\begin{array}{l}\text{PSI}\\ g_{\mathbf{(2,1)}}\end{array}$}
				\TR{$\begin{array}{l}\text{CP}\\ g_{\mathbf{(2,2)}}\end{array}$}
				\TR{$\begin{array}{l}\text{R}\\ g_{\mathbf{(2,3)}}\end{array}$}
				}				
				\pstree{ \TR{$\begin{array}{l}\textbf{En}\\ G_{\mathbf{3}}\end{array}$}}{
				\TR{$\begin{array}{l}\text{TPR}\\ g_{\mathbf{(3,1)}}\end{array}$}
				\TR{$\begin{array}{l}\text{IWU}\\ g_{\mathbf{(3,2)}}\end{array}$}
				\TR{$\begin{array}{l}\text{LS}\\ g_{\mathbf{(3,3)}}\end{array}$}
					\TR{$\begin{array}{l}\text{OU}\\ g_{\mathbf{(3,4)}}\end{array}$}
				}		
				}}
\end{figure}

The elementary subcriteria are described in Table~\ref{criteriadefinition}, with their different scale of evaluation, and a description of whether they have to be maximized (denoted by $\uparrow$) or minimized (denoted by $\downarrow$). \\
The whole hierarchy of criteria is displayed in Figure~\ref{hierarchy_sicily}.

\begin{table}[!htb]
\caption{Description of the elementary subcriteria}\label{criteriadefinition}\smallskip
\centering
\begin{small}
\resizebox{1\textwidth}{!}{
\begin{tabular}{llcc} 
\hline
\hline
\textbf{Elementary subriterion} & \textbf{Description}  & \textbf{Unit}  & \textbf{Preference}\\
\hline
\hline
Owned Homes/Total Homes (OH) & ratio between the number of total homes owned by the residents and the total number of  homes & $[0,1]$ & $\downarrow$\\[0,7mm]
Agricultural Age Structure (AA)  &ratio between farms conducted by young farmers and the total number of farms & $[0,1]$  & $\uparrow$\\[0,7mm]
Bed Spaces in Rural Tourism Accommodation/km$^2$ (BS) & number of bed spaces in rural tourism accommodation per unit of area & $1/km^2$  & $\uparrow$\\[0,7mm]
\hline
Population Scattering Index (PSI)&  ratio between the population living in individual homes, villages and small towns, and the total population &  $[0,1]$ & $\uparrow$\\[0,7mm]
Concentrated Population/km$^2$ (CP)  &  concentrated population per unit of area &  $1/km^2$  &$\downarrow$ \\[0,7mm]
Total Residents (R) & total number of residents &   &$\downarrow$ \\[0,7mm]
\hline
Typical Products Rate (TPR) & ratio between the number of typical products from agriculture in the municipality and the total number of typical products of Sicily & $[0,1]$   & $\uparrow$\\
Irrigation Water Use/UAA (IWU) &  amount of water used for irrigation per unit (Hectare) of Utilized Agricultural Area (UAA)&   $m^3/$Hectare & $\downarrow$ \\
Livestock Standard units/UAA (LS) & number of animals breeding converted into Livestock Standard Units per unit of UAA &   LSU/UUA & $\downarrow$\\
Organic UAA/UAA (OU) &ratio between organic UAA and total UAA & $[0,1]$  & $\uparrow$\\								
\hline
\hline
\end{tabular}
}
\end{small}
\end{table}

Now assume that the DM -- e.g., a policy maker who has a preference structure oriented towards rural sustainability -- provides the following preference information on the considered elementary criteria and on the macro-criteria (notice that each statement is translated into a constraint, which is written in round brackets): 

\begin{description}
\item[(i)] \textbf{So} is more important than \textbf{En}, which, in turn, is more important than \textbf{Ec} $\big(\varphi_{\mathbf{0}}(\text{So})\geq \varphi_{\mathbf{0}}(\text{En})+\varepsilon$ and $\varphi_{\mathbf{0}}(\text{En})\geq \varphi_{\mathbf{0}}(\text{Ec})+\varepsilon$\big);
\item[(ii)] at a comprehensive level, AA is more important than TRP $\left(\varphi_{\mathbf{0}}(\text{AA})\geq \varphi_{\mathbf{0}}(\text{TRP})+\varepsilon \right)$;
\item[(iii)]  with respect to \textbf{Ec}, OU is more important than LS $\left(\varphi_{\mathbf{1}}^2(\text{OU})\geq \varphi_{\mathbf{1}}^2(\text{LS})+\varepsilon \right)$;
\item[(iv)] with respect to \textbf{So}, PSI is more important than CP $\left(\varphi_{\mathbf{2}}^2(\text{PSI})\geq \varphi_{\mathbf{2}}^2(\text{CP})+\varepsilon \right)$;
\item[(v)] at a comprehensive level, CP is more important than OH $\left(\varphi_{\mathbf{0}}^2(\text{CP})\geq \varphi_{\mathbf{0}}^2(\text{OH})+\varepsilon \right)$;
\item[(vi)] with respect to \textbf{En}, LS and OU are positively interacting $\left(\varphi_{\mathbf{3}}^2(\text{LS}, \text{OU})\geq \varepsilon \right)$;
\item[(vii)] BS and TPR are positively interacting $\left(\varphi_{\mathbf{0}}^2(\text{BS},\text{TPR})\geq \varepsilon \right)$;
\item[(viii)] OH and CP are negatively interacting $\left(\varphi_{\mathbf{0}}^2(\text{OH}, \text{CP})\leq -\varepsilon \right)$.
\end{description}

Before implementing the proposed methodology, we perform a normalization procedure, because applications of the Choquet integral always require that all evaluations are expressed on the same scale. 
Specifically, if criterion $g_i$ has a decreasing direction of preference, then the following expression is used:
$$
\overline{g_i}(a)=\frac{\max_i-g_i(a)}{\max_i-\min_i}, 
$$
where $g_i(a) $ is the evaluation of $a$ on criterion $g_i$, $\max_i$ and $\min_i$ are the maximum and the minimum evaluations of alternatives on $g_i$, and $\overline{g_i}(a)$ is the normalized value. 
On the other hand, whenever criterion $g_i$  has an increasing direction of preference, the following dual formula is adopted:

$$
\overline{g_i}(a)=\frac{g_i(a)-\min_i}{\max_i-\min_i}. 
$$ 

\noindent Notice that the normalization chosen for this application is the most used among the different methods available in the literature on sustainability~\cite{Singh:2012}.

Regrettably, the application of NAROR at the comprehensive level and at each macro-criterion (social, economic and environmental) yields results that are insufficient to get a good insight into our decision problem. Indeed, beyond the relation of weak dominance, the application of NAROR only adds few new pairs to the necessary preference relations at the comprehensive level, as well as on the macro-criteria.  Specifically, at the comprehensive level, NAROR provides the following additional necessary relations: 
$$
a_{35} \!\succsim^N \!\!a_{45}\:, \;\; a_{44} \!\succsim^N \!\!a_{41}\:, \;\;a_{9} \!\succsim^N \!\!a_{39}\:,\;\; a_{32} \!\succsim^N \!\!a_{39} \:,\;\; a_{51} \!\succsim^N \!\! a_{39} \:,\;\; a_{21} \! \succsim^N \!\! a_{39} \:,\;\; a_{36} \!\succsim^N \!\! a_{39} \:,\;\; a_{36} \!\succsim^N \!\! a_{15}\:.
$$
Further, with respect to the macro-criterion \textbf{En}, NAROR only provides the additional necessary preference relation $a_{11} \! \succsim_{\mathbf{3}}^{N} \! a_{15}$.

\medskip

Thus, in order to get a more refined understanding of the case study under examination, we implement SMAA. To start, we compute the best and the worst position reachable from each municipality: to that end, we consider the whole set of capacities compatible with the preferences provided by the DM, as well as the three rank positions presenting the highest rank acceptability indices (hence the three most frequent positions taken by that municipality). 

In order to compare municipalities, we shall give much attention to the reference ranking on \textbf{So}, since this macro-criterion is the most important according to the DM's evaluation (cf.\ preference information (i) provided by the DM). Specifically, we rank the municipalities with respect to their highest rank acceptability index placing. Therefore, for instance, in the first positions one finds those municipalities that present the highest rank acceptability index for the first positions: see Table~\ref{RAI_So}, where we report the three highest rank acceptabilities indices, as well as the best and worst positions of the first seven and the last three municipalities in the considered ranking. For the sake of completeness, in Table~\ref{RAI_GL}-\ref{RAI_Ec}-\ref{RAI_En}, the same computations are performed at a comprehensive level, on macro-criterion \textbf{Ec}, and on macro-criterion \textbf{En}, respectively. 

\begin{table}[!h]
\begin{center}
\caption{The three highest rank acceptability indices and the best and worst positions for the municipalities ranked (a) at a comprehensive level, (b) according to macro-criterion \textbf{Ec}, (c) according to macro-criterion \textbf{So}, and (d) according to macro-criterion \textbf{En}.\label{RAI}}
\subtable[Comprehensive level\label{RAI_GL}]{
\resizebox{0.7\textwidth}{!}{
\begin{tabular}{ccccccccccc}
    \hline
		\hline
\rule[-2mm]{0mm}{0.8cm}		
Municipality &  $high_1$  & $b_{\mathbf{0}}^{high_1}(\cdot)\left(\%\right)$ & $high_2$ & $b_{\mathbf{0}}^{high_2}(\cdot)\left(\% \right)$ & $high_3$ & $b_{\mathbf{0}}^{high_3}(\cdot)\left( \% \right)$ &Best & $b_{\mathbf{0}}^{Best}(\cdot)\left( \% \right)$ & Worst & $b_{\mathbf{0}}^{Worst}(\cdot)\left( \%\right)$ \\[1ex]
\hline
\hline
$a_{25}$ &    2     & 46.44 & 1     & 32.41 & 3     & 14.36 & 1     & 32.41 & 13    & 0.00 \\[1ex] 
$a_{20}$ &     1     & 58.01 & 2     & 36.58 & 3     & 4.83  & 1     & 58.01 & 7     & 0.00 \\[1ex]  
$a_{12}$ &     4     & 34.63 & 5     & 25.32 & 3     & 20.46 & 2     & 0.37  & 14    & 0.00 \\[1ex]  
$a_{26}$ &     12    & 18.04 & 11    & 17.47 & 10    & 14.48 & 4     & 0.07  & 27    & 0.00 \\[1ex]  
$a_{34}$ &    8     & 18.70 & 10    & 18.19 & 9     & 16.31 & 3     & 0.01  & 20    & 0.00 \\[1ex]  
$a_{44}$ &   5     & 18.13 & 4     & 15.22 & 7     & 13.35 & 1     & 0.05  & 27    & 0.00 \\[1ex]  
$a_{13}$ &   13    & 21.42 & 14    & 18.46 & 15    & 14.52 & 5     & 0.00  & 32    & 0.00 \\[1ex]  
$\cdots$ &           &     &                 &  &                &                & & & &  $\cdots$         \\[1ex]   
$a_{4}$    &  49    & 25.26 & 48    & 23.54 & 50    & 15.20 & 42    & 0.00  & 51    & 6.24 \\[1ex]  
$a_{47}$   &  46    & 25.26 & 47    & 19.25 & 48    & 15.90 & 29    & 0.00  & 51    & 0.10 \\[1ex]  
$a_{2}$    &  51    & 76.33 & 50    & 19.32 & 49    & 3.93  & 46    & 0.00  & 51    & 76.33 \\[1ex]  
    \hline
		\hline
\end{tabular}  
}}
\subtable[Economic (Ec)\label{RAI_Ec}]{%
\resizebox{0.7\textwidth}{!}{
\begin{tabular}{ccccccccccc}
    \hline
		\hline
\rule[-2mm]{0mm}{0.8cm}		
Municipality &  $high_1$  & $b_{\mathbf{1}}^{high_1}(\cdot)\left(\%\right)$ & $high_2$ & $b_{\mathbf{1}}^{high_2}(\cdot)\left(\% \right)$ & $high_3$ & $b_{\mathbf{1}}^{high_3}(\cdot)\left( \% \right)$ &Best & $b_{\mathbf{1}}^{Best}(\cdot)\left( \% \right)$ & Worst & $b_{\mathbf{1}}^{Worst}(\cdot)\left( \%\right)$ \\[1ex]
\hline
\hline
 $a_{25}$ &   27    & 17.17 & 28    & 14.96 & 26    & 14.62 & 18    & 0.01  & 43    & 0.00 \\[1ex]
$a_{20}$ &    2     & 26.11 & 3     & 18.08 & 4     & 17.60 & 1     & 14.76 & 18    & 0.00 \\[1ex]
$a_{12}$ &    9     & 27.81 & 10    & 27.41 & 8     & 18.82 & 6     & 2.19  & 25    & 0.00 \\[1ex]
 $a_{26}$ &   28    & 11.01 & 27    & 8.70  & 29    & 8.26  & 7     & 0.00  & 45    & 0.01 \\[1ex]
 $a_{34}$ &   36    & 26.39 & 37    & 23.81 & 33    & 9.64  & 22    & 0.00  & 43    & 0.00 \\[1ex]
  $a_{44}$ & 13    & 8.60  & 6     & 6.60  & 14    & 5.05  & 3     & 1.26  & 48    & 0.00 \\[1ex]
  $a_{13}$ & 46    & 27.72 & 47    & 25.83 & 48    & 24.41 & 23    & 0.00  & 49    & 3.30 \\[1ex]
$\cdots$ &           &     &                 &  &                &                & & & &  $\cdots$         \\[1ex]   
  $a_{4}$ &  31    & 12.48 & 33    & 10.86 & 36    & 10.54 & 6     & 0.01  & 41    & 0.01 \\[1ex]
  $a_{47}$ &  33    & 13.26 & 31    & 13.25 & 29    & 10.53 & 19    & 0.00  & 47    & 0.00 \\[1ex]
 $a_{2}$ &   43    & 21.13 & 40    & 16.45 & 39    & 14.19 & 33    & 0.00  & 49    & 0.02 \\[1ex]
    \hline
		\hline
\end{tabular}
}}
\subtable[Social (So)\label{RAI_So}]{%
\resizebox{0.7\textwidth}{!}{
\begin{tabular}{ccccccccccc}
    \hline
		\hline
\rule[-2mm]{0mm}{0.7cm}		
Municipality &  $high_1$  & $b_{\mathbf{2}}^{high_1}(\cdot)\left(\%\right)$ & $high_2$ & $b_{\mathbf{2}}^{high_2}(\cdot)\left(\% \right)$ & $high_3$ & $b_{\mathbf{2}}^{high_3}(\cdot)\left( \% \right)$ &Best & $b_{\mathbf{2}}^{Best}(\cdot)\left( \% \right)$ & Worst & $b_{\mathbf{2}}^{Worst}(\cdot)\left( \%\right)$ \\[1ex]
\hline
\hline
   $a_{25}$ &  1     & 100.00 &     &   &     &   & 1     & 100.00 & 1     & 100.00 \\[1ex]
   $a_{20}$ & 2     & 99.55 & 3     & 0.38  & 5     & 0.06  & 2     & 99.55 & 8     & 0.00 \\[1ex]
  $a_{12}$ &  3     & 96.40 & 4     & 2.69  & 5     & 0.43  & 3     & 96.40 & 13    & 0.00 \\[1ex]
   $a_{26}$ & 4     & 92.87 & 6     & 3.87  & 5     & 2.30  & 2     & 0.00  & 15    & 0.00 \\[1ex]
  $a_{34}$ &  5     & 94.76 & 6     & 3.86  & 4     & 1.18  & 4     & 1.18  & 9     & 0.00 \\[1ex]
   $a_{44}$ &  6     & 91.52 & 3     & 3.12  & 4     & 2.88  & 2     & 0.45  & 7     & 0.00 \\[1ex]
   $a_{13}$ & 7     & 79.62 & 8     & 18.21 & 9     & 0.74  & 2     & 0.00  & 11    & 0.01 \\[1ex]
$\cdots$ &           &     &                 &  &                &                & & & &  $\cdots$ \\[1ex]
  $a_{4}$ &  49    & 36.25 & 51    & 30.30 & 47    & 14.17 & 43    & 0.00  & 51    & 30.30 \\[1ex]
  $a_{47}$ &50    & 49.54 & 49    & 29.39 & 48    & 7.91  & 41    & 0.00  & 51    & 6.03 \\[1ex]
  $a_{2}$ & 51    & 63.64 & 50    & 31.44 & 49    & 4.45  & 46    & 0.02  & 51    & 63.64 \\[1ex]
    \hline
		\hline
\end{tabular} 
 }}
\subtable[Environmental (En)\label{RAI_En}]{%
\resizebox{0.7\textwidth}{!}{
\begin{tabular}{ccccccccccc}
    \hline
		\hline
\rule[-2mm]{0mm}{0.7cm}		
Municipality &  $high_1$  & $b_{\mathbf{3}}^{high_1}(\cdot)\left(\%\right)$ & $high_2$ & $b_{\mathbf{3}}^{high_2}(\cdot)\left(\% \right)$ & $high_3$ & $b_{\mathbf{3}}^{high_3}(\cdot)\left( \% \right)$ &Best & $b_{\mathbf{3}}^{Best}(\cdot)\left( \% \right)$ & Worst & $b_{\mathbf{3}}^{Worst}(\cdot)\left( \%\right)$ \\[1ex]
\hline
\hline
$a_{25}$ &    17    & 29.95 & 18    & 29.50 & 19    & 12.15 & 14    & 0.55  & 30    & 0.00 \\[1ex] 
$a_{20}$ &     14    & 9.95  & 4     & 9.16  & 13    & 7.98  & 4     & 9.16  & 33    & 0.06 \\[1ex]
$a_{12}$ &    15    & 10.53 & 16    & 8.78  & 14    & 6.95  & 5     & 0.28  & 34    & 0.12 \\[1ex] 
$a_{26}$ &    26    & 22.17 & 25    & 21.84 & 27    & 19.68 & 15    & 0.00  & 31    & 0.00 \\[1ex] 
$a_{34}$ &    8     & 33.61 & 9     & 28.38 & 7     & 14.36 & 4     & 0.00  & 18    & 0.01 \\[1ex] 
  $a_{44}$ &  12    & 17.83 & 6     & 15.35 & 11    & 13.15 & 3     & 0.00  & 18    & 0.09 \\[1ex] 
     $a_{13}$ & 9     & 28.79 & 10    & 20.82 & 8     & 17.62 & 3     & 0.00  & 20    & 0.00 \\[1ex]
$\cdots$ &           &     &                 &  &                &                & & & &  $\cdots$         \\[1ex]   
 $a_{4}$ &   42    & 45.58 & 41    & 13.16 & 43    & 12.04 & 22    & 0.00  & 47    & 0.00 \\[1ex] 
 $a_{47}$ &   26    & 10.46 & 27    & 9.40  & 37    & 7.93  & 9     & 0.00  & 50    & 0.00 \\[1ex]
 $a_{2}$ &    44    & 30.55 & 43    & 30.28 & 45    & 12.84 & 18    & 0.00  & 50    & 0.04 \\[1ex]
 \hline
 \hline
 \end{tabular}
}}
\end{center}
\end{table}

Let us first analyze Table~\ref{RAI_So}, which, according to the importance given by the DM to macro-criterion \textbf{So} may be considered the most relevant of all four. 
As a matter of fact, we employ the results obtained for \textbf{So} as a reference for the analysis of all other aspects. 
More specifically, we shall restrict our analysis to (i) the top seven alternatives w.r.t.\ \textbf{So}, and (ii) the last three alternatives w.r.t.\ \textbf{So}.

\begin{itemize}
\item[(i)] $a_{25}$ (Milo) is surely the best municipality,  since its rank acceptability index for the first position is equal to $100\%$. From a social point of view, also $a_{20}$ (Maniace), $a_{34}$ (Ragalna) and $a_{44}$ (Sant'Alfio) appear quite good, since they have high rank acceptability indices for the first rank positions, and take always a position in the intervals $[2,8]$, $[4,9]$, and $[2,7]$, respectively (see Table~\ref{RAI_So}). At the same time, even if $a_{12}$ (Castel di Iudica) and $a_{26}$ (Mineo) have high rank acceptabilities for the third and fourth places, respectively, their rank position varies in the intervals $[3,13]$ for Castel di Iudica, and $[2,15]$ for Mineo. 
\item[(ii)] The last three municipalities present the highest rank acceptability indices for the lowest positions. In particular, $a_4$ (Acireale) is very frequently in the 49th position; it is the last ranked municipality in almost one third of the cases, and its best position is the 43rd. Furthermore, $a_{47}$ (Tremestieri Etneo) presents the highest rank acceptability index for the last but one place in the ranking (49.54\%), and it may reach at its best the 41st position, even if with a negligible frequency. Finally, $a_2$ (Aci Catena) is the municipality presenting the highest rank acceptability for the last position ($63.64\%$), and its second and third highest rank acceptability indices are for the the 50th and the 49th positions. The best position it reaches is the 46th, even if the rank acceptability index for this position is basically zero (0.02$\%$). 
\end{itemize}

In what follows, we provide some comments on the rankings obtained in the other three cases, and compare  the relative performances with those displayed in the reference rankings on macro-criterion \textbf{So}. 
We start with the rankings given in Table~\ref{RAI_Ec}, computed according to macro-criterion \textbf{Ec}.

\begin{itemize}
\item[(i)] Differently from what happens with respect to macro-criterion \textbf{So}, $a_{25}$ (Milo) appears to be quite unstable with respect to macro-criterion \textbf{Ec}: indeed, its highest rank acceptability index, corresponding to the $27$th position, is very low ($17.17\%$), and the interval of variation for its positions is quite large, ranging from the 18th to the 43rd place. Instead, $a_{20}$ (Maniace) is quite good on economic aspects, since its three highest rank acceptability indices are those corresponding to the positions going from the 2nd to the 4th and it is always placed between the $1$st and the $18$th position. Furthermore, $a_{12}$ (Castel di Iudica) is evaluated quite well on \textbf{Ec}, since it presents its three highest rank acceptability indices for the positions going from the $8$th  to the $10$th; however, its range of positions is rather large, being within the interval $[6,25]$. Finally, it is interesting that municipalities $a_{34}$ (Raddusa) and $a_{13}$ (Fiumefreddo di Sicilia), which showed good performances on \textbf{So}, always have their three highest rank acceptability indices for positions at the bottom of the ranking with respect to macro-criterion \textbf{Ec}. 
\item[(ii)]  The last three municipalities on \textbf{So} -- namely, $a_4$, $a_{47}$, $a_{2}$ -- improve their positions with respect to macrocriterion \textbf{Ec}, even if they are still placed at the bottom of the ranking. In particular, $a_4$ and $a_{47}$ can reach a position in the first half of the ranking, but this happens with very marginal frequencies. Moreover, $a_4$ has a frequency of the $0.01\%$ of being in the 6th position, whereas $a_{47}$ has the 19th as its best rank-position, even if with an almost null frequency.
\end{itemize}

Analyzing the same municipalities with respect to macro-criterion \textbf{En} (see Table~\ref{RAI_En}), the following observations can be done.

\begin{itemize}
\item[(i)] $a_{34}$ is the best among the first seven considered municipalities, insofar as it has the highest rank acceptability indices for the positions going from the 7th to the 9th. Moreover, this municipality can  fill the fourth place in the ranking even if with a negligible frequency, and its last position is the 18th again with a frequency very close to zero. The other two municipalities presenting quite stable results are $a_{13}$ and $a_{44}$, since they have the highest rank acceptability for positions not lower than 12th, and they can never reach a position lower than the 20th. 
\item[(ii)] Among the three last municipalities -- namely, $a_{4}$, $a_{47}$, $a_{2}$ -- only $a_{47}$ improves its position, since its most frequent position is the 26th even if with a low rank acceptability index ($10.46\%$). Moreover, even if $a_{4}$ and $a_{2}$ can reach the 22nd and the 18th as their highest positions, they are very frequently in the last part of the ranking. Indeed, $a_4$ presents its highest rank acceptability indices for positions 42nd, 41st and 43rd, whereas $a_2$ has the 44th, 43rd and 45th as the positions for which it has the highest rank acceptability indices.
\end{itemize}

Finally, examining Table~\ref{RAI_GL}, the following remarks can be done at a comprehensive level.

\begin{itemize}
\item[(i)] $a_{20}$ and $a_{25}$ can be considered the best municipalities, since they have the three highest rank acceptability indices for the first three positions. Moreover, $a_{20}$ can be placed at the top of this ranking, because its worst rank-position is the 7th against the $13$th of $a_{25}$. Furthermore, $a_{12}$ can be ranked just after the first two, since its most frequent positions are the 4th, the 5th and the 3rd with frequencies $34.63\%$, $25.32\%$ and $20.46\%$, respectively. 
\item[(ii)] The three last municipalities do not improve the position they had in the \textbf{So} ranking. Indeed, all of them present their highest rank acceptability indices for the lowest positions in the ranking. Moreover, $a_2$ takes the last position ($51$st) with the first highest rank acceptability equal to $76.33\%$, and the other two municipalities can reach the last position, too (in particular, $a_4$ with a frequency of $6.24\%$). 
\end{itemize}

In order to get a ranking of the considered municipalities with respect to \textbf{Ec}, \textbf{So}, \textbf{En} and at the comprehensive level, we have calculated the \textit{barycenter} of the M\"{o}bius representation of capacities compatible with the preferences provided by the DM. 
Their values are shown in Table ~\ref{Average_coeff}.

From this table, one can notice that, without considering interaction between indicators, the first three more important elementary criteria are the ones relative to the macro-criterion \textbf{So}, i.e., PSI, R and CP, with $m(\left\{PSI\right\})> m(\left\{R\right\}) >m(\left\{CP\right\})$. Some positive and negative interactions between criteria -- apart from the preference information explicitly given by the DM -- can be deduced from the same table: for instance, the negative and positive interaction, respectively, between PSI and R on one hand, and LS and OU on the other hand).

\begin{table}[!h]
\caption{Barycenter  of the M\"{o}bius representations of compatible capacities}\label{Average_coeff}\smallskip
\centering
\begin{small}
\resizebox{1\textwidth}{!}{
\begin{tabular}{ccccccccccc} 
\hline
\hline
\rule[-2mm]{0mm}{0.7cm}
 $m(\left\{OH\right\})$ & $m(\left\{AA\right\})$ & $m(\left\{BS\right\})$ & $m(\left\{PSI\right\})$ & $m(\left\{CP\right\})$ & $m(\left\{R\right\})$ & $m(\left\{TPR\right\})$ & $m(\left\{IWU\right\})$ & $m(\left\{LS\right\})$ &  $m(\left\{OU\right\})$ & $m(\left\{OH,AA\right\})$ \\
    0.0871 & 0.1103 & 0.0869 & 0.1508 & 0.1187 & 0.1218 & 0.0695 & 0.0882 & 0.0649 & 0.0948 & -0.0001\\
\hline
\hline
\rule[-2mm]{0mm}{0.7cm}
 $m(\left\{OH,BS\right\})$ & $m(\left\{OH,PSI\right\})$ & $m(\left\{OH,CP\right\})$  &   $m(\left\{OH,R\right\})$  &  $m(\left\{OH,TPR\right\})$  &  $m(\left\{OH,IWU\right\})$  & $m(\left\{OH,LS\right\})$  &  $m(\left\{OH,OU\right\})$  &  $m(\left\{AA,BS\right\})$  & $m(\left\{AA,PSI\right\})$  &$m(\left\{AA,CP\right\})$ \\
    -0.0003 & -0.0002 & -0.0163 & -0.0001 & 0.0004 & -0.0002 & 0.0000 & 0.0006 & -0.0006	 &-0.0044	& 0.0006\\
\hline
\hline
\rule[-2mm]{0mm}{0.7cm}
  $m(\left\{AA,R\right\})$ & $m(\left\{AA,PR\right\})$ & $m(\left\{AA,IWU\right\})$ & $m(\left\{AA,LS\right\})$ & $m(\left\{AA,OU\right\})$ & $m(\left\{BS,PSI\right\})$ & $m(\left\{BS,CP\right\})$ & $m(\left\{BS,R\right\})$ & $m(\left\{BS,TPR\right\})$ & $m(\left\{BS,IWU\right\})$ & $m(\left\{BS,LS\right\})$ \\
    0.0005 & -0.0005 & -0.0007 & 0.0000 & 0.0005 & 0.0000	& 0.0004	& -0.0014	& 0.0133	& -0.0005	&-0.0003\\
\hline
\hline
\rule[-2mm]{0mm}{0.7cm}
$m(\left\{BS,OU\right\})$ & $m(\left\{PSI,CP\right\})$ & $m(\left\{PSI,R\right\})$ & $m(\left\{PSI,TPR\right\})$ & $m(\left\{PSI,IWU\right\})$ & $m(\left\{PSI,LS\right\})$ & $m(\left\{PSI,OU\right\})$ & $m(\left\{CP,R\right\})$ & $m(\left\{CP,TPR\right\})$ & $m(\left\{CP,IWU\right\})$ & $m(\left\{CP,LS\right\})$  \\
-0.0001  & 0.0002	& -0.0008 &	-0.0010 &	0.0011	& 0.0005 &	0.0010 & -0.0001 &	0.0003 &	-0.0002 &	0.0001 \\
\hline
\hline
\rule[-2mm]{0mm}{0.7cm}
$m(\left\{CP,OU\right\})$ & $m(\left\{R,TPR\right\})$ & $m(\left\{R,IWU\right\})$ & $m(\left\{R,LS\right\})$ & $m(\left\{R,OU\right\})$ & $m(\left\{TPR,IWU\right\})$ & $m(\left\{TPR,LS\right\})$ & $m(\left\{TPR,OU\right\})$ & $m(\left\{IWU,LS\right\})$ & $m(\left\{IWU,OU\right\})$ & $m(\left\{LS,OU\right\})$  \\
0.0003 &  -0.0002	& 0.0004	& -0.0002 &	-0.0003 & 0.0001 &	0.0000	&-0.0002 & 0.0004	& 0.0006 & 0.0144\\
\hline
\hline
\end{tabular}
}
\end{small}
\end{table}

Computing the Choquet integral value for each municipality with the barycenter of the M\"{o}bius representation of all compatible capacities, we obtain four complete rankings of municipalities at the comprehensive level as well as at the macro-criteria level. In Tables ~\ref{Rank0}, ~\ref{Rank1},~\ref{Rank2} and~\ref{Rank3}, we show the places taken by the municipalities under consideration  by using the capacity in terms of  M\"{o}bius representation displayed in Table \ref{Average_coeff}. 

\begin{table}[!h]
\begin{center}
\caption{Positions taken by the municipalities in the rankings obtained by using the  M\"{o}bius representation shown in Table \ref{Average_coeff} \label{Rank}}
\subtable[Comprehensive level\label{Rank0}]{%
\resizebox{0.45\textwidth}{!}{
\begin{tabular}{cc}
\hline
\hline
\rule[-2mm]{0mm}{0.7cm}
Municipality & Position in the complete ranking    \\
\hline
\hline
Milo	 & $2$\\
Maniace	& $1$\\
Castel di Iudica	& $4$\\
Mineo	& $12$\\
Ragalna	& $10$\\
Sant'Alfio	& $5$\\
Castiglione di Sicilia	& $13$\\
$\cdots$  &  $\cdots$  \\
Acireale	& $49$ \\
Tremestieri Etneo	& $46$\\
Aci Catena	& $51$\\
\hline
\hline
\end{tabular}
}
}
\subtable[Economic (\textbf{Ec}) \label{Rank1}]{%
\resizebox{0.45\textwidth}{!}{
\begin{tabular}{cc}
\hline
\hline
\rule[-2mm]{0mm}{0.7cm}
Municipality & Position in the complete ranking    \\
\hline
\hline
Milo	 & $28$\\
Maniace	& $3$\\
Castel di Iudica	& $8$\\
Mineo	& $27$\\
Ragalna	& $36$\\
Sant'Alfio	& $14$\\
Castiglione di Sicilia	& $13$\\
$\cdots$  &  $\cdots$  \\
Acireale	& $49$ \\
Tremestieri Etneo	& $46$\\
Aci Catena	& $51$\\
\hline
\hline
\end{tabular}
}
}
\subtable[Social (\textbf{So})\label{Rank2}]{%
\resizebox{0.45\textwidth}{!}{
\begin{tabular}{cc}
    \hline
		\hline
\rule[-2mm]{0mm}{0.7cm}		
Municipality & Position in the complete ranking    \\
\hline
\hline
Milo	 & $1$\\
Maniace	& $2$\\
Castel di Iudica	& $3$\\
Mineo	& $4$\\
Ragalna	& $5$\\
Sant'Alfio	& $6$\\
Castiglione di Sicilia	& $7$\\
$\cdots$  &  $\cdots$  \\
Acireale	& $49$\\
Tremestieri Etneo	& $50$\\
Aci Catena	& $51$\\		
\hline
\hline
\end{tabular}
}
}
\subtable[Environmental (\textbf{En})\label{Rank3}]{%
\resizebox{0.45\textwidth}{!}{
\begin{tabular}{cc}
    \hline
		\hline
\rule[-2mm]{0mm}{0.7cm}		
Municipality & Position in the complete ranking    \\
\hline
\hline
Milo	 & $18$\\
Maniace	& $13$\\
Castel di Iudica	& $16$\\
Mineo	& $24$\\
Ragalna	& $8$\\
Sant'Alfio	& $11$\\
Castiglione di Sicilia	& $9$\\
$\cdots$  &  $\cdots$  \\
Acireale	& $42$\\
Tremestieri Etneo	& $33$\\
Aci Catena	& $43$\\			
\hline
\hline
\end{tabular}
}
}
\end{center}
\end{table}

As expected, the position taken by a municipality  depends on the criterion we are considering. For example, with respect to macro-criteria \textbf{Ec} and \textbf{En}, Mineo takes positions 27 and 24, respectively, whereas its position becomes the $12$th at a comprehensive level. Similarly, Milo is respectively ranked $28$th and $18$th on macro-criteria \textbf{Ec} and \textbf{En}; instead, at a comprehensive level it takes the $2$nd position.  

The obtained results highlight the features of the employed approach, which provides insights into evaluation of sustainable development from two different perspectives: for each node of the hierarchy of the criteria, and at a global level. The interested reader can download the file containing complete results by clicking on the following link: \href{www.antoniocorrente.it/salvo/images/DataFIR3.xls}{data}.

\section{Conclusions}\label{conclusion}
We proposed a methodology to construct composite indices taking into account the following relevant points:

\begin{itemize}
\item interaction of elementary indicators;
\item hierarchical organization of the considered indicators;
\item participation of experts, stakeholders and policy makers to the construction of the composite indices;
\item consideration of robustness concerns related to the stability of the results with respect to the variability of the weights assigned to considered indicators.
\end{itemize}

We suggested to handle these points simultaneously by applying a multiple criteria decision analysis methodology, that is the hierarchical-SMAA-Choquet integral, that has been recently proposed in literature. This approach:
\begin{itemize}
\item aggregates the indicators by means of the  Choquet integral preference model that takes into account the possible interactions between them,
\item applies the Robust Ordinal Regression (ROR) and the Stochastic Multicriteria Acceptability Analysis (SMAA) to handle robustness concerns related to the consideration of the whole set of weights assigned to composite indicators and compatible with preferences expressed by experts, stakeholders and policy makers,
\item adopts Multiple Criteria Hierarchy Process (MCHP) to aggregate indcators not only at the comprehensive level, but also at all the intermediate levels of the hierarchy.
\end{itemize}

We have discussed how this approach is well suited for composite indices related to sustainable development assessment and we have shown how it works on a case study related to rural sustainable development. 

We believe that the methodology we presented has interesting properties that suggest its use in all the domains in which composite indices are adopted or can be adopted, with a specific relevance to sustainable development. Consequently we plan to apply this methodology to construct composite indices in different areas of interest.   

From a more theoretical point of view, we envisage to apply other MCDA methodologies to the same case study. In fact, a comparison of the obtained results would highlight the different features of the employed methods.

\section*{Acknowledgments}
The authors wish to acknowledge the funding by the ``FIR of the University of Catania BCAEA3, New developments in Multiple Criteria Decision Aiding (MCDA) and their application to territorial competitiveness" and of the research project ``Data analytics for entrepreneurial ecosystems, sustainable development and wellbeing indices" of the Department of Economics and Business of the University of Catania. Salvatore Greco has also benefited of the fund ``Chance" of the University of Catania, while Silvia Angilella, Alfio Giarlotta and Marcella Rizzo have also benefited of the fund ``FFABR" of the Ministry of Education, University and Research of the Italian government.

\section*{References}
\newcommand{\noopsort}[1]{} \newcommand{\printfirst}[2]{#1}
  \newcommand{\singleletter}[1]{#1} \newcommand{\switchargs}[2]{#2#1}

\end{document}